\documentclass[journal]{IEEEtran}
\usepackage{xspace}        
\usepackage{mathtools}
\usepackage[hidelinks]{hyperref}
\usepackage{bbm}
\usepackage{url}
\usepackage{cleveref}
\usepackage{dsfont} 
\usepackage{siunitx}
\usepackage{booktabs}       
\usepackage{multirow}
\usepackage{algorithmic}
\usepackage[ruled,vlined]{algorithm2e}       
\usepackage{amsthm,amsmath,amssymb}
\usepackage{color}
\usepackage{comment}
\usepackage{tikz}
\usepackage{pgfplots}
\usepackage{bbm}
\usepackage{cite}
\usepackage{graphicx}
\usepackage{textcomp}
\usepackage{epstopdf}
\usepackage{bm}
\usepackage{color}
\usepackage{mathrsfs}
\usepackage{tabu}

\newtheorem{lemma}{Lemma}
\newtheorem{theorem}{Theorem}
\newtheorem{assumption}{Assumption}

\ifCLASSINFOpdf
\else
\fi
\begin{document}
%
\title{ Asymmetric Learning in Convex Games }
%
%
%
\author{Zifan Wang, Xinlei Yi, Yi Shen, Michael M. Zavlanos, and Karl H. Johansson
\thanks{* This work was supported in part by Swedish Research Council Distinguished Professor Grant 2017-01078, Knut and Alice Wallenberg Foundation, Wallenberg Scholar Grant, the Swedish Strategic Research Foundation SUCCESS Grant, AFOSR under award \#FA9550-19-1-0169, and  NSF under award CNS-1932011. (Corresponding author: Xinlei Yi.)}
\thanks{Zifan Wang and Karl H. Johansson are with the Division of Decision and Control Systems, School of Electrical Enginnering and Computer Science, KTH Royal Institute of Technology, and also with Digital Futures, SE-10044 Stockholm, Sweden. Email: \{zifanw,kallej\}@kth.se.}
\thanks{Xinlei Yi is with Shanghai Institute of Intelligent Science and Technology, Tongji University, Shanghai, 201804, China, and is also with Lab for Information \& Decision Systems, Massachusetts Institute of Technology, Cambridge, MA 02139, USA. Email: xinleiyi@tongji.edu.cn.}
\thanks{Yi Shen and Michael M. Zavlanos are with the Department of Mechanical Engineering and Materials Science, Duke University, Durham, NC, USA. Email: \{yi.shen478, michael.zavlanos\}@duke.edu}
}

%
%

\allowdisplaybreaks

\markboth{Journal of \LaTeX\ Class Files,~Vol.~14, No.~8, August~2015}%
{Shell \MakeLowercase{\textit{et al.}}: Asymmetric Learning in Convex Games}
%



\maketitle

\begin{abstract}
This paper considers convex games involving multiple agents that aim to minimize their own cost functions using locally available information. A common assumption in the study of such games is that the agents are symmetric, meaning that they have access to the same type of information. Here we lift this assumption, which is often violated in practice, and instead consider asymmetric agents; specifically, we assume some agents have access to first-order gradient information and others have access to the zeroth-order oracles (cost function evaluations). We propose an asymmetric learning algorithm that combines the agent information mechanisms. We analyze the regret and Nash equilibrium convergence of this algorithm for convex and strongly monotone games, respectively. Specifically, we show that our algorithm always performs between pure first- and zeroth-order methods, and can match the performance of these two extremes by adjusting the number of agents with access to zeroth-order oracles. Therefore, our algorithm incorporates the pure first- and zeroth-order methods as special cases. We provide numerical experiments on a market problem for both deterministic and risk-averse games to demonstrate the performance of the proposed algorithm.
\end{abstract}


%
\IEEEpeerreviewmaketitle

\section{Introduction}

Convex optimization \cite{boyd2004convex,zinkevich2003online,hazan2016introduction,nesterov2018lectures} is widely applicable in various fields such as economics, engineering, and machine learning.
Recently, convex optimization has been employed in multi-agent games with applications in traffic routing \cite{sessa2019no} and market optimization \cite{wang2022risk}.
In these applications, agents are usually assumed rational with the goal to minimize their own cost functions by leveraging limited information received from the environment, which falls into the category of convex games \cite{zhang2021multi,bravo2018bandit}.

The performance of optimization algorithms for convex games is typically evaluated using the notion of regret \cite{hazan2016introduction}, which captures the difference between  agents' sequential actions and the  achievable  best actions in hindsight. 
An algorithm is said to achieve no-regret learning if the regret of a sequence of actions generated by the algorithm is sub-linear in the total number of episodes $T$, meaning that agents are able to eventually learn the best actions. 
Another important measure of performance in convex games is that of a  Nash equilibrium, which  is defined as a  point at which no agent has incentive to change its decision.
Recently, there is a growing literature in game theory focusing on designing algorithms that achieve 
no-regret learning \cite{mertikopoulos2018cycles,wang2022zeroth,sessa2019no,hsieh2021adaptive} or Nash equilibrium convergence \cite{belgioioso2020distributed,lee2021last,bervoets2020learning,bravo2018bandit,tatarenko2020bandit,tatarenko2018learning,paccagnan2018nash,frihauf2011nash,tembine2013risk}. 
Common in these works is the assumption that the agents have access to  similar or symmetric type of information.
However, in many real-world settings, agents are asymmetric, meaning that they have access to different types or amounts of information; 
for instance, in financial markets, investment banks have access to more information compared to 
individual investors. 
Moreover, in Cournot competition, where multiple companies aim to maximize their own profit, certain dominant companies have the capability to observe their competitors' strategies while others do not have this advantage. 
Information asymmetry is also evident in security systems \cite{etesami2019dynamic,xu2015exploring}, where an attacker operates in obscurity and has the capability to observe the actions of the defender. Conversely, the defender operates with transparency, often without access to the attacker's information or actions.
In these cases, the asymmetric information  gives rise to asymmetric update strategies for the agents.

To the best of our knowledge, asymmetric learning in convex games has not been explored in the literature.
Most closely related to our study is  symmetric learning in convex games \cite{bravo2018bandit,hsieh2022no,lin2020finite,golowich2020tight,mertikopoulos2020gradient,mertikopoulos2021equilibrium,mertikopoulos2019learning,narang2022learning,duvocelle2023multiagent,drusvyatskiy2022improved,balduzzi2018mechanics,chasnov2020convergence,mazumdar2020gradient}, where authors
have proposed methods for no-regret learning and/or Nash equilibrium convergence. 
Specifically, when gradient information is available,  \cite{hsieh2022no} develops optimistic gradient methods for continuous games with noisy gradient estimates that achieve constant regret under multiplicative noise.
Similarly, \cite{lin2020finite} proposes a first-order gradient descent algorithm for $\lambda$-coercive games with unconstrained continuous action sets, which attains the last-iterate convergence to a Nash equilibrium. 
Assuming bandit feedback, i.e., agents only have access to zeroth-order oracles, \cite{bravo2018bandit} shows Nash equilibrium convergence  for  strongly monotone games, which is a special class of convex games. The convergence rate of the zeroth-order method in \cite{bravo2018bandit} is further improved in \cite{drusvyatskiy2022improved} relying on the additional assumption that the Jacobian of the gradient function is Lipschitz continuous.
Common in these works is that the agents perform symmetric updates using the same kind of information. The methods cannot be directly extended to asymmetric agents.

There are a few works that analyze asymmetric learning in other classes of games with different types of asymmetry \cite{kononen2004asymmetric,tampubolon2020information,xu2015exploring,li2022commitment,arslantas2023convergence}.
For example, \cite{kononen2004asymmetric} considers information asymmetry in Stacklberg games, where one agent (leader) can observe the action of the other agent (follower), and proposes a learning method that utilizes the theory of Markov games. 
Subsequently, \cite{tampubolon2020information} proposes an asymmetric Q-learning algorithm for two-agent Markov games and discusses the existence of Nash equilibria and convergence. 
In \cite{li2022commitment}, Bayesian Stackelberg games are analyzed under double-sided
information asymmetry where the leader hides its action from the follower and the follower holds information about its payoff private. 
It is shown that the leader can improve its payoff by strategically revealing part of the action to the follower.
We note that all the works above focus on Stackelberg games or Markov games with two agents acting \textit{consecutively}, which differ from the convex games considered by us.

In this paper, we consider asymmetry in the form of gradient  available to the agents during the learning process. Specifically, we assume that some agents have access to zeroth-order oracles, while other agents additionally have access to first-order gradient information. 
This situation may arise when some agents can observe other agents' actions and thus can compute the first-order gradient while the remaining agents have only access to cost function evaluations at each episode.
Asymmetry also arises in stochastic games involving both risk-neutral and risk-averse agents. When  Conditional Value at Risk (CVaR) is used as a risk measure, the CVaR gradient can rarely be explicitly derived even if the  
form of the stochastic cost function is known \cite{rockafellar2000optimization}. The risk-averse agents cannot easily have access to gradient information unlike the risk-neutral agents.

Our main contributions are summarized as follows.
\begin{enumerate}
\item 
We develop a novel framework that relies on asymmetric information to learn optimal actions in convex games. In this framework, agents belong to two distinct groups: agents that only have access to zeroth-order oracles and update their actions using zeroth-order optimization techniques, and agents that have access to first-order gradient information and update their actions accordingly.
\item 
While the asymmetric setting complicates the system dynamics, we theoretically 
show that no-regret learning is achieved for every agent for convex games, and last-iterate Nash equilibrium convergence is guaranteed for strongly monotone games.
We show that the performance of the proposed asymmetric learning algorithm lies between the pure first- and zeroth-order methods.
\item 
Experimentally, we validate our algorithm on online markets, specifically, a deterministic and a risk-averse Cournot game. In the latter case, the cost of each agent is stochastic and agents may be risk-neutral or risk-averse to avoid catastrophically high costs. We show that when there are risk-neutral agents, our asymmetric learning algorithm converges faster compared to the pure zeroth-order method in \cite{wang2022risk}.
\end{enumerate}

To the best of our knowledge, this is the first paper to address asymmetric agent updates in convex games with theoretical results on regret analysis and Nash equilibrium convergence. 
Perhaps closest to the method proposed here is the one in  \cite{niu2022dish}, which addresses a distributed consensus optimization problem where the agents  perform either Newton- or gradient-type updates. Linear convergence is shown for strongly convex objective functions regardless of updates.
We note that the games considered here are different than the consensus optimization problems in \cite{niu2022dish}, as is also the type of asymmetric information.
As a result, the techniques in \cite{niu2022dish} cannot be applied to analyze the asymmetric games considered here.

The rest of this paper is organized as follows. In Section~\ref{sec:1_problem}, we formally define asymmetric games and provide some basic notation. 
The asymmetric learning algorithm is proposed in Section~\ref{sec:algorithm}. 
In Section~\ref{sec:convex},  we analyze regret of the algorithm for convex games.
Section~\ref{sec:strongly} provides Nash equilibrium convergence analysis for strongly monotone games of the asymmetric algorithm.
Section~\ref{sec:experiments} experimentally illustrates the algorithm in the application to deterministic and risk-averse  Cournot games. Finally, we conclude the paper in Section~\ref{sec:conclusion}.
The detailed proofs of the results presented in this paper can be found in the arXiv paper \cite{wang2023asymmetric}.

\section{Problem Definition}
\label{sec:1_problem}

Consider a repeated game with $N$ non-cooperative agents, whose goals are to learn the best actions that minimize their own cost functions. 
Let $\mathcal{N} = \{1, 2, \ldots, N\}$ denote the set of agent indices.
For each agent  $i\in \mathcal{N}$, the cost function is $\mathcal{C}_i(x_i,x_{-i}) : \mathcal{X} \rightarrow \mathbb{R}$, where $x_i \in \mathcal{X}_i$ is the action of agent $i$ and  $x_{-i}$  the actions of all agents except for agent $i$. The joint action space is defined as $\mathcal{X} =\Pi_{i=1}^N\mathcal{X}_i$, where $ \mathcal{X}_i \subset \mathbb{R}^d$, $d>0$, is a convex set.
We write $x:=(x_1,\ldots,x_N)$ as the  collection of all agents' actions. 
Throughout this paper, we use the notation $(x,y)$ to denote the concatenated vector $[x^{\rm{T}},y^{\rm{T}}]^{\rm{T}}$.

The goal of each agent $i$ is to minimize its individual cost function, i.e.,
\begin{align}\label{eq:def:game}
    \mathop{{\rm{minimize}}}_{x_i \in \mathcal{X}_i} \; \mathcal{C}_i(x_i,x_{-i}).
\end{align}
The game \eqref{eq:def:game} is defined as a convex game when each agent's cost function is convex in its individual action.
We consider the following class of convex games.
\begin{assumption}\label{assump:convex}
For each agent $i\in \mathcal{N}$,
$\mathcal{C}_i(x_i,x_{-i})$ is convex in $x_i$ for all $x_{-i} \in \mathcal{X}_{-i}$, where $\mathcal{X}_{-i} =\Pi_{j=1,j\neq i }^N\mathcal{X}_j$. Moreover, we assume that $\mathcal{X}_i$ contains the ball with radius $r$ centered at the origin and has a bounded diameter $D>0$, for all $i\in \mathcal{N}$. 
\end{assumption} 
Assumption~\ref{assump:convex} is common in the literature, see, e.g., \cite{drusvyatskiy2022improved}.
As shown in \cite{rosen1965existence}, there always exists at least one Nash equilibrium in a convex game \eqref{eq:def:game}. We denote by $x^{*}$ such a Nash equilibrium. By definition, we have that 
\begin{align*}
    \qquad \qquad \mathcal{C}_i(x^{*})\leq \mathcal{C}_i(x_i,x_{-i}^{*}), \;  {\rm{for \; all}} \; x_i \in \mathcal{X}_i, i\in \mathcal{N}.
\end{align*} 
At a Nash equilibrium point, no agent can reduce its own cost by unilaterally changing its individual action.
Since every agent's cost function is convex in its own action, the Nash equilibrium can be characterized using the first-order optimality condition 
\begin{align*}
    \qquad \langle \nabla_i \mathcal{C}_i(x^{*}), x_i - x_i^{*} \rangle \geq 0, \;   {\rm{for \; all}}\;  x_i \in \mathcal{X}_i, \; i\in \mathcal{N},
\end{align*}
where we write $\nabla_{i} \mathcal{C}_i(x)$ instead of $\nabla_{x_i} \mathcal{C}_i(x)$ for the ease of notation and the symbol $\nabla_{x_i}$ means taking the gradient with respect to $x_i$.

We consider the case when agents have access to asymmetric information. 
Specifically, there are $N_z \geq 0$ agents that only have access to zeroth-order oracles, and $N-N_z$ agents that have access to first-order gradient information. 
We index the agents that only have access to zeroth-order oracles by $i_z$, where $i_z \in \mathcal{N}_z:=\{ 1,\ldots,N_z\}$, and the agents that have access to first-order gradient information by $i_f$, where $i_f\in \mathcal{N}_f:=\{ N_z+1,\ldots,N\}$. 
When $N_z = 0$, we denote $\mathcal{N}_z = \emptyset$, and when $N_z = N$, $\mathcal{N}_f = \emptyset$.
Let the vectors $x_z := (x_1, \ldots, x_{N_z})$ and $x_f := (x_{N_z + 1}, \ldots, x_N)$ represent the action profiles of agents that have access to zeroth-order oracles and first-order gradient information, respectively.
Finally, we use the notation $-i_z$ to represent all agents that have access to zeroth-order oracles except for the agent $i_z$.

We make the following assumptions on the game \eqref{eq:def:game}.
\begin{assumption}\label{assump:U}
For each agent $i\in \mathcal{N}$ and $x\in \mathcal{X}$, there exists $U>0$ such that $|\mathcal{C}_i(x)|\leq U$.
\end{assumption}

\begin{assumption}\label{assump:L0}
For each agent $i\in \mathcal{N}$, there exists $L_0>0$ such that  $\mathcal{C}_i(x)$ is $L_0$-Lipschitz continuous in $x$, i.e., 
$$|\mathcal{C}_i(x) - \mathcal{C}_i(y)| \leq L_0 \left\| x-y \right\| , \; \forall x,y\in \mathcal{X}.$$
\end{assumption}
\begin{assumption}\label{assump:L1}
For each agent $i\in \mathcal{N}$, there exists $L_1>0$ such that $ \nabla_i \mathcal{C}_i(x)$ is $L_1$-Lipschitz continuous in $x$, i.e.,
$$
\left\|\nabla_i \mathcal{C}_i(x) - \nabla_i \mathcal{C}_i(y)\right\| \leq L_1 \left\| x-y\right\|,  \; \forall x,y\in \mathcal{X}.
$$
\end{assumption}
These assumptions are common in the literature and hold in many applications, e.g., Cournot games and Kelly Auctions \cite{bravo2018bandit,lin2021doubly,duvocelle2023multiagent}.

Convergence analysis for games with multiple Nash equilibria is in general hard. For this reason, recent research has focused on the so-called strongly monotone games, which are shown in \cite{rosen1965existence} to admit a unique Nash equilibrium.
The game \eqref{eq:def:game} is said to be $m$-strongly monotone with $m>0$, if for all $ x,x'\in \mathcal{X}$ it satisfies
\begin{align}\label{eq:strong_monotone}
    \sum_{i=1}^N \langle \nabla_i \mathcal{C}_i(x) -\nabla_i \mathcal{C}_i(x'),x_i-x_i' \rangle \geq m \left\|x -x' \right\|^2.
\end{align}

For convex games, the ability of agents to efficiently learn their optimal actions can be quantified using the notion of regret, which captures the difference between the agents' sequential actions and the  achievable  best actions in hindsight. 
Given the sequences of actions $\{x_{i,t} \}_{t=1}^T, i=1,\ldots,N$, the regret of agent~$i$ is  defined as 
\begin{align}\label{eq:def:regret:game}
    {\rm{R}}_i(T)= \sum_{t=1}^T \mathcal{C}_i(x_{i,t},x_{-i,t}) - \mathop{\rm{min}}_{x_i \in \mathcal{X}_i} \sum_{t=1}^T\mathcal{C}_i(x_i,x_{-i,t}).
\end{align}
An algorithm is said to be no-regret if the regret of each agent is sub-linear in $T$.

Our goal is to design an asymmetric learning algorithm to solve the game \eqref{eq:def:game}, when the zeroth-order agents $i_z \in \mathcal{N}_z$ and the first-order agents $i_f \in \mathcal{N}_f$  update their actions using the specific type information available to them. We aim to show that the proposed algorithm achieves  no-regret learning and Nash equilibrium convergence for convex  and strongly monotone games, respectively.

\section{An Asymmetric Learning Algorithm}\label{sec:algorithm}
\begin{algorithm}[tb]
\caption{Asymmetric learning} \label{alg:algorithm_asym}
\begin{algorithmic}[1]
    \STATE \textbf{Input}:  Initial value $x_1$, positive sequences $\eta_{f,t}$, $\eta_{z,t}$, $\delta_t$, parameters $N_z$, $N$, $T$.
    \FOR{$ {\text{episode}} \; t=1,\ldots,T$}
        \STATE Agents play their actions:
        \FOR{agent $i=1,\ldots,N$}
            \STATE Agent $i_z \in \mathcal{N}_z$ samples $u_{i_z,t} \in \mathbb{S}$ and plays $\hat{x}_{i_z,t} = x_{i_z,t}+\delta_t u_{i_z,t}$ 
            \STATE Agent $i_f \in \mathcal{N}_f$ plays $ x_{i_f,t}$
        \ENDFOR
        \STATE Agents receive information:
        \FOR{agent $i=1,\ldots,N$}
            \STATE Agent $i_z \in \mathcal{N}_z$ obtains $\mathcal{C}_{i_z}(\hat{x}_{z,t},x_{f,t})$
            \STATE Agent $i_f \in \mathcal{N}_f$ obtains $\nabla_{i_f} \mathcal{C}_{i_f}(x_{f,t},\hat{x}_{z,t})$
        \ENDFOR
        \STATE Agents perform updates:
        \FOR{agent $i=1,\ldots,N$}
            \STATE Agent $i_z \in \mathcal{N}_z$ updates according to \eqref{eq:asy:update:ZO}
            \STATE Agent $i_f \in \mathcal{N}_f$ updates according to \eqref{eq:asy:update:FO}
        \ENDFOR
    \ENDFOR
    \STATE  \textbf{Output}: $\{x_{t}\}$.
\end{algorithmic}
\end{algorithm}
In this section, we propose an asymmetric learning algorithm that allows the agents to update their actions using first-order or zeroth-order gradient information, whichever is available to them. 

By saying that agent $i$ has access to the first-order gradient information, we mean that agent $i$ can compute the gradient $\nabla_i \mathcal{C}_i$. If agent $i$ only has access to zeroth-order oracle, then the available information for agent $i$ is its function evaluation $\mathcal{C}_i$.
We assume that agents are rational and always use first-order gradient information to update their actions when available, while those with only zeroth-order gradient information use this information for updating.

The detailed algorithm is presented as Algorithm~\ref{alg:algorithm_asym}. 
At each episode $t$, each agent $i_z \in \mathcal{N}_z $ that has access to a zeroth-order oracle perturbs its action by an amount of $\delta_t u_{i_z,t}$ and plays the perturbed action $\hat{x}_{i_z,t} = x_{i_z,t}+\delta_t u_{i_z,t}$, where $u_{i_z,t} \in \mathbb{S}$ is a random perturbation direction sampled from a unit sphere $\mathbb{S}\subset \mathbb{R}^{d}$ and $\delta_t$ is the size of this perturbation.  
In contrast, each agent $i_f \in \mathcal{N}_f$ with first-order gradient information plays the action $x_{i_f,t}$. Here we collect the actions of all agents that have access to zeroth-order oracles in a vector $\hat{x}_{z,t}:= (\hat{x}_{1,t}, \ldots, \hat{x}_{N_z,t})$, and the actions of all agents that have access to first-order gradient information in a vector $x_{f,t}:= (x_{N_z+1,t},\ldots,x_{N,t})$.
After all agents have played their actions, they receive their own asymmetric information. 
Specifically, given the played action profile $(\hat{x}_{z,t},x_{f,t})$ at time step $t$, each agent $i_z$ obtains the returned function evaluation $\mathcal{C}_{i_z}(\hat{x}_{z,t},x_{f,t})$, while each agent $i_f$ gets the first-order gradient information $\nabla_{i_f} \mathcal{C}_{i_f}(x_{f,t},\hat{x}_{z,t})$. Then, each agent $i_z$ performs the following projected update
\begin{align}\label{eq:asy:update:ZO}
    x_{i_z,t+1} &\leftarrow \mathcal{P}_{\mathcal{X}_{i_z}^{\delta_t}} \Big( x_{i_z,t} - \eta_{z,t} {g}_{i_z,t}\Big), 
\end{align}
where the step size $\eta_{z,t}$ of agents with zeroth-order oracles will be designed later.
The zeroth-order gradient estimate is constructed as
\begin{align}
    g_{i_z,t}=\frac{d}{\delta_t} \mathcal{C}_{i_z}(\hat{x}_{z,t},x_{f,t}) u_{i_z,t}.
\end{align} The projection set is defined as
$ \mathcal{X}_{i_z}^{\delta_t} = \{ x_{i_z}\in \mathcal{X}_{i_z} \vert \frac{1}{1-\delta_t/r}x_{i_z} \in \mathcal{X}_{i_z}\}$. The projection step guarantees the feasibility of the sampled action $\hat{x}_{i_z,t}$, since 
\begin{align*}
    &\bigg(1-\frac{\delta_t}{r}\bigg) \mathcal{X}_{i_z} \oplus \delta_t \mathbb{S} \\
    &= \bigg(1-\frac{\delta_t}{r}\bigg) \mathcal{X}_{i_z} \oplus \frac{\delta_t}{r} r   \mathbb{S}\\
    & \subseteq \bigg(1-\frac{\delta_t}{r}\bigg)\mathcal{X}_{i_z} \oplus \frac{\delta_t}{r} \mathcal{X}_{i_z} = \mathcal{X}_{i_z}.
\end{align*}
Here, $\oplus$ denotes the Minkowski sum of two sets. 
Each agent $i_f$ performs the action update
\begin{align}\label{eq:asy:update:FO}
    x_{i_f,t+1} &\leftarrow \mathcal{P}_{\mathcal{X}_{i_f}} \Big( x_{i_f,t} - \eta_{f,t} \nabla_{i_f} \mathcal{C}_{i_f}(x_{f,t},\hat{x}_{z,t}) \Big),
\end{align}
where $\eta_{f,t}$ is the step size of agents with first-order gradient information and will be designed later.

In our framework, each agent $i_z\in \mathcal{N}_z$ with no access to first-order gradient information plays a randomly perturbed action $\hat{x}_{i_z,t}$, as in standard zeroth-order optimization \cite{flaxman2004online}. Then, the agents with access to the first-order gradients obtain gradient information at the random point $(x_{f,t},\hat{x}_{z,t})$ rather than $(x_{f,t},{x}_{z,t})$. 
In other words, agents with access to zeroth-order oracles add randomness to the action updates of agents with first-order gradient information. This randomness does not exist in the pure first-order method where all agents obtain gradient information at a deterministic point $x_t$.
Consequently, such asymmetric updates lead to more complex agent interactions.

Although Algorithm 1 may appear to be a combination of the pure first- and zeroth-order methods, the analysis of these two methods cannot be directly applied to our algorithm. On one hand, as stated above, the presence of agents with zeroth-order oracles introduces extra randomness to the action updates of agents with first-order gradient information, thereby increasing the whole system's complexity. 
On the other hand, agents with first-order oracles introduce discoordination for the whole system. Specifically, the dynamics of first-order agents directly affect the cost function values, thereby affecting the gradient estimates for agents with zeroth-order oracles.
Consequently, the constructed gradient estimates become biased of the smoothed functions in pure zeroth-order case, rendering traditional zeroth-order techniques inapplicable here.

\section{Convergence Results for Convex Games}\label{sec:convex}
In this section, we analyze the convergence of Algorithm~\ref{alg:algorithm_asym} given that the game \eqref{eq:def:game} is convex.

The standard smoothed functions commonly used in analyzing pure zeroth-order methods rely on the premise that all agents obtain gradient estimates based on zeroth-order optimization techniques. However, the presence of agents with first-order gradient information in our asymmetric game renders the premise invalid, making these smooth functions unsuitable for our analysis.
To solve this problem, we construct new smoothed functions tailored for such asymmetric games.
Specifically,  we define
$ \mathcal{C}_{i_z,t}(x_z):= \mathcal{C}_{i_z}(x_z,x_{f,t})$ and the smoothed function $ \mathcal{C}^{\delta_t}_{i_z,t}(x_z):= \mathop{\mathbb{E}}_{w_{i_z} \sim \mathbb{B},u_{-i_z} \sim \mathbb{S}_{-i_z}} \mathcal{C}_{i_z,t}(x_{i_z}+ \delta_t w_{i_z}, x_{-i_z} + \delta_t u_{-i_z})$ which serves as an approximation of $\mathcal{C}_{i_z,t}(x_z)$. Here, $\mathbb{B}$, $\mathbb{S}$ denote the unit ball and unit sphere in $\mathbb{R}^{d}$, respectively, and $\mathbb{S}_{-{i_z}}:=\Pi_{j=1}^{N_z-1} \mathbb{S}$.
The smoothed cost function in symmetric games \cite{bravo2018bandit} is smoothed over actions of all agents, while our smoothed cost function is only smoothed over actions of zeroth-order agents. This modification lays the foundation for our subsequent analysis. When all agents use zeroth-order optimization, the smoothed function reduces to the one in \cite{bravo2018bandit}.
It can be shown that the function $\mathcal{C}^{\delta_t}_{i_z,t}(x_z)$ has the following properties \cite{wang2023asymmetric}. 
%
\begin{lemma}\label{lemma:func:property}
Let Assumptions \ref{assump:convex} and \ref{assump:L0} hold. Then, we have that for all $t\geq 1$
\begin{enumerate}
    \item $\mathcal{C}_{i_z,t}^{\delta_t}(x_{z})$ is convex in $x_{i_z}$;
    \item $\mathcal{C}_{i_z,t}^{\delta_t}(x_z)$ is $L_0$-Lipschitz continuous in $x_z$;
    \item $|\mathcal{C}_{i_z,t}^{\delta_t}(x_z)-\mathcal{C}_{i_z,t}(x_z)|\leq \delta_t L_0 \sqrt{N_z}$;
    \item $\mathbb{E}\left[ \frac{d}{\delta_t} \mathcal{C}_{i_z,t}(\hat{x}_{z}) u_{i_z} \right] = \nabla_{i_z} \mathcal{C}^{\delta_t}_{i_z,t}(x_z)$.
\end{enumerate}
\end{lemma}
The smoothed function $\mathcal{C}_{i_z,t}^{\delta_t}$ is defined as a time-varying cost function of the actions of all agents that have access to zeroth-order oracles. Despite its time-varying nature, Lemma~\ref{lemma:func:property} shows that some properties of the function $\mathcal{C}_{i_z,t}^{\delta_t}$ still hold for any $x_{f,t}$.
Note that the last property in Lemma~\ref{lemma:func:property} shows that the term $\frac{d}{\delta_t} \mathcal{C}_{i_z,t}(\hat{x}_{z}) u_{i_z}$ is an unbiased estimate of the gradient of the smoothed function $\mathcal{C}^{\delta_t}_{i_z,t}(x_z)$.

We use the notion of regret to measure the performance of our algorithm. Given a sequence of agent actions $\{\hat{x}_{i_z,t} \}_{t=1}^T$, $i_z \in \mathcal{N}_z$, and $\{x_{i_f,t} \}_{t=1}^T$, $i_f \in \mathcal{N}_f$ generated by Algorithm~\ref{alg:algorithm_asym}, 
\eqref{eq:def:regret:game} yields the regret of agent $i_z \in \mathcal{N}_z$ as
\begin{align*}
    {\rm{R}}_{i_z}(T)= \mathbb{E} \Big[ &\sum_{t=1}^T \mathcal{C}_{i_z}(\hat{x}_{z,t},x_{f,t}) \nonumber \\
    &- \mathop{\rm{min}}_{x_{i_z} \in \mathcal{X}_{i_z}} \sum_{t=1}^T \mathcal{C}_{i_z}(x_{i_z}, \hat{x}_{-i_z,t},x_{f,t}) \Big],
\end{align*}
and agent $i_f \in \mathcal{N}_f$ as 
\begin{align*}
    {\rm{R}}_{i_f}(T)= & \mathbb{E} \Big[ \sum_{t=1}^T \mathcal{C}_{i_f}(x_{f,t},\hat{x}_{z,t}) \nonumber \\
    &- \mathop{\rm{min}}_{x_{i_f} \in \mathcal{X}_{i_f}} \sum_{t=1}^T \mathcal{C}_{i_f}(x_{i_f}, x_{-i_f,t},\hat{x}_{z,t}) \Big].
\end{align*}

We note that the definition of regret depends on actions of other agents, and the asymmetric updates of other agents make regret analysis more complex. 
By appropriately selecting the parameters $\eta_{f,t}$, $\eta_{z,t}$ and $\delta_t$, we show that Algorithm~\ref{alg:algorithm_asym} achieves sub-linear regret. The formal result is presented in the following theorem, in which the notion $\mathcal{O}$ hides all constant factors except for $N$, $N_z$ and $T$. 
Due to space limitations, the detailed proof is provided in 
\cite{wang2023asymmetric}. 
\begin{theorem}\label{theorem:No-regret}
Let Assumptions~\ref{assump:convex}--\ref{assump:L0} hold, and select  $\delta_t =  \frac{1}{N_z^{1/4} t^{1/4}} $, $\eta_{z,t} = \frac{1}{N_z^{1/4} t^{3/4}}$, $\eta_{f,t} = \frac{1}{ \sqrt{t}}$, for all $ t\geq 1$. Then, for any $T \geq 1$, Algorithm~\ref{alg:algorithm_asym} achieves regrets
\begin{align}
    &{\rm{R}}_{i_z}(T) = \mathcal{O}\left( N_z^{1/4}T^{3/4} \right), \; \forall i_z \in \mathcal{N}_z   ; \nonumber \\
    &{\rm{R}}_{i_f}(T) = \mathcal{O}\left( T^{1/2} \right), \; \forall i_f \in \mathcal{N}_f.
\end{align}
\end{theorem}

The design of the positive sequences $\eta_{z,t}$, $\eta_{f,t}$ and $\delta_{t}$ does not require the knowledge of the total number of episodes $T$. Using these sequences,
Theorem~\ref{theorem:No-regret} shows that Algorithm~\ref{alg:algorithm_asym} achieves sub-linear regret for all agents, although agents that have access to first-order gradient information have a smaller regret bound in terms of  $T$ compared to agents accessing zeroth-order oracles.

Using similar techniques as in \cite{hazan2016introduction}, one can show that the pure first-order method achieves a regret of $\mathcal{O}\left( T^{1/2} \right)$  while the pure zeroth-order method achieves a regret of $ \mathcal{O}\left( N^{1/4}T^{3/4} \right)$. 
We observe that the regret for agents with first-order gradient information is not affected by the asymmetric setting. This is due to the facts that other agents' varying actions can be deemed as a source of varying cost functions and first-order gradient descent update is no-regret under varying cost functions. However, the agents with zeroth-order oracles are affected by this discoordination induced by the asymmetric setting. Specifically, the dynamics of the first-order agents directly affect the cost function values, thereby affecting the zeroth-order gradient estimates. Despite this, we show that the constructed gradient estimate is an unbiased estimate of the gradient of the newly constructed smooth function, and the regret bound of zeroth-order agents achieved by Algorithm~\ref{alg:algorithm_asym} is always better than that achieved by the pure zeroth-order method.
Therefore, Algorithm~\ref{alg:algorithm_asym} not only ensures that the agents accessing first-order gradient information inherit the lower regret bound $\mathcal{O}\left( T^{1/2} \right)$ of the pure first-order method, but also improves the regret bound of the agents accessing zeroth-order oracles, compared to the regret achieved by  the pure zeroth-order method. 
Let the total regret denote the sum of every agent's regret, which we denote by $R = \sum_{i_z \in \mathcal{N}_z} R_{i_z } +\sum_{i_f \in \mathcal{N}_f} R_{i_f }$. The total regret achieved by Algorithm~\ref{alg:algorithm_asym} is lower bounded by that of the pure zeroth-order method and upper bounded by that of the pure first-order method.

\section{Convergence Results for Strongly Monotone Games}\label{sec:strongly}

In this section, we analyze Nash equilibrium convergence given that the game \eqref{eq:def:game} is strongly monotone. It is well-known  that the Nash equilibrium is unique in strongly monotone games \cite{rosen1965existence}, which we denote by $x^{*}$. In what follows, we provide the last-iterate Nash equilibrium convergence result for Algorithm~\ref{alg:algorithm_asym}.

\begin{theorem}\label{theorem:NE}
Suppose that the game~\eqref{eq:def:game} is $m$-strongly monotone. Let Assumptions~\ref{assump:convex}--\ref{assump:L1} hold, and select $\delta_t = \frac{N_z^{1/6}}{N^{1/3} t^{1/3}} $ , $\eta_{z,t} = \eta_{f,t} =\frac{1}{mt}$, for all $ t= 1,\ldots,T$.  Then, Algorithm~\ref{alg:algorithm_asym} satisfies
\begin{align}\label{eq:NE:convergence:asym}
    &\mathbb{E}\left\| x_T - x^{*} \right\|^2 \nonumber \\
    &= \mathcal{O}\Big( N_z^{2/3} N^{2/3} T^{-1/3} +  (N-N_z)T^{-1}\Big).
\end{align}
\end{theorem}

\noindent \emph{Proof Sketch:}
We analyze the convergence to the Nash equilibrium by separately examining the error dynamics for first- and zeroth-order agents and then combining these analyses together.
Specifically, for first-order agents, we study the evolution of the squared error $\left\|x_{i_f,t+1} - x_{i_f}^*\right\|^2$ leveraging the corresponding update rule \eqref{eq:asy:update:FO}. 
For zeroth-order agents, due to the varying parameters of $\delta_t$ in the smoothed function for zeroth-order agents, we first analyze the evolution of $\left\| x_{i_z,t+1} - x_{i_z,t}^*\right\|$ leveraging the update rule \eqref{eq:asy:update:ZO}, where $x_{i_z,t}^* = (1-\frac{\delta_t}{r}) x_{i_z}^*$. Then, we analyze $\left\| x_{i_z,t+1} - x_{i_z}^*\right\|$ by establishing a bound on $\left\| x_{i_z}^* - x_{i_z,t}^*\right\|^2$.

Although the dynamics of the two agent types are interconnected, we can bound the influence of one agent type on the other by leveraging the properties established in Lemma~\ref{lemma:func:property}.
Next, we combine the dynamics of all agents to analyze the behavior at the group level, i.e., the convergence toward the Nash equilibrium.
By applying the strong monotonicity condition, we characterize the evolution of the Nash equilibrium error, quantified by $\left\|x_t - x^*\right\|^2$.
Finally, we prove convergence using induction.
$\hfill \blacksquare $

Theorem~\ref{theorem:NE} shows that Algorithm~\ref{alg:algorithm_asym} achieves the last-iterate Nash equilibrium convergence for strongly monotone games with diminishing smoothing parameters $\delta_t$ and diminishing step sizes $\eta_{f,t}$ and $\eta_{z,t}$. Note that the design of diminishing parameters does not require the information of the total number of episodes $T$. Moreover, under Assumption~\ref{assump:L0}, Nash equilibrium convergence implies sub-linear regret for each agent since the regret is upper bounded by the sum of differences between their actions and Nash equilibrium.

When $N_z = 0$, meaning that all agents use the first-order gradient-descent update, the convergence rate $O (N T^{-1})$ matches the result of learning in strongly monotone games with stochastic gradient information \cite{jordan2024adaptive}.
In the other extreme case $N_z=N$, i.e., all agents have only access to zeroth-order oracles, the term $(N-N_z) T^{-1}$ disappears and the convergence rate $O(N^{4/3} T^{-1/3})$ matches the result of zeroth-order learning in games \cite{bravo2018bandit}.
We observe that the convergence rate of Algorithm~\ref{alg:algorithm_asym} lies between those of first- and zeroth-order methods.
It is worth noting that the convergence rate of $O(T^{-1/3})$ for zeroth-order optimization in games is suboptimal, as indicated in \cite{lin2021doubly}. Therefore, there is potential to improve the convergence rate in \eqref{eq:NE:convergence:asym}. However, achieving this in the asymmetric setting is challenging, and we leave it for future work.

To end this section, we would like to briefly explain the challenges in the proof of Theorem~\ref{theorem:NE}. Algorithm~\ref{alg:algorithm_asym} is a combination of the pure first- and zeroth-order methods. However, the Nash equilibrium convergence analysis of Algorithm~\ref{alg:algorithm_asym} cannot be obtained by directly combining the results of these two pure  methods, since the dynamics of the two groups of agents, i.e., agents with first-order gradient information and zeroth-order oracles, are coupled together. On the one hand, agents with access to zeroth-order oracles spread their randomness to  the dynamics of the agents with first-order gradient information.
On the other hand, the presence of first-order agents slightly disrupt the coordination among those agents with zeroth-order oracles aiming to minimize the smoothed cost function.
This coupled dynamics, combined with the setting of diminishing parameters, complicates the analysis of Nash equilibrium convergence of the whole system.

\section{Numerical Experiments}\label{sec:experiments}
In this section we use a Cournot game to illustrate the performance of Algorithm~\ref{alg:algorithm_asym}. Consider a market problem with $N$ agents. Suppose that each agent $i$ supplies the market with quantity $x_i$ and the total supply of all agents decides the price of the goods in the market. Given the price, each agent $i$ has the cost $\mathcal{C}_i(x)$. Each agent aims to minimize its own cost through repeated learning in the game.
In what follows, we consider two cases, a deterministic game and a stochastic game with risk-averse agents.

\subsection{Deterministic Cournot Game}
Consider a market with $N=10$ agents. The cost function of each agent is given by 
$$\mathcal{C}_i(x) = x_i(\frac{a_i x_i}{2} + b_i \sum_{j \neq i}x_{j} - e_i)+ 1,$$ where $a_i>0$, $b_i$, $e_i$ are constant parameters. It is easy to verify that $$\nabla_i \mathcal{C}_i(x) = a_i x_i + b_i x_{-i} -e_i.$$ 
The parameters are selected as 
\begin{align*}
    a &= [2, 2, 1.5, 1.8, 2,1.8, 2, 1.4,1.8, 2], \nonumber \\
    b &= [0.2, 0.3, 0.3, 0.2, 0.3, 0.2, 0.3, 0.2, 0.3, 0.3], \nonumber \\
    e &= [1.8, 1.9, 1.5, 1.6, 1.8, 1.3, 1.2, 1.5, 1.8, 1.6].
\end{align*}
It can be verified that the game is $m$-strongly monotone with $m=1.284$, and the Nash equilibrium is $$x^{*} = [0.57,0.44,0.29,0.52,0.38,0.33,0.026,0.61,0.43,0.26].$$
The projection set is defined as  $\mathcal{X}_i = [0,3]$.
\begin{figure}[tb]
\begin{center}
\centerline{\includegraphics[width=1\columnwidth]{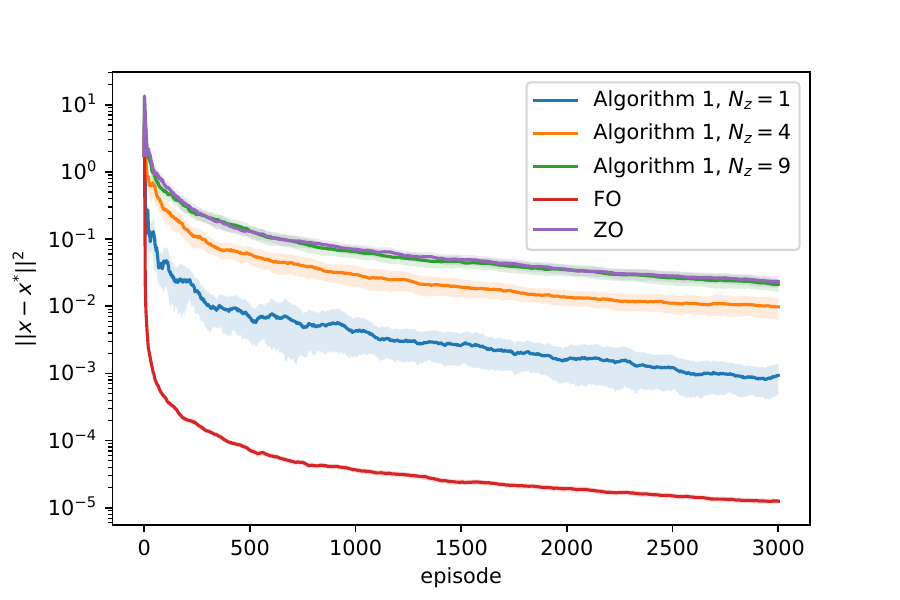}}
\caption{Nash equilibrium convergence achieved by  the pure first-order method (FO) and the pure zeroth-order method (ZO), and our asymmetric learning algorithm with different values of $N_z$. The solid lines and shades are averages and standard deviations over 50 runs.}
\label{fig_NE}
\end{center}

\end{figure}

\begin{figure}[tb]
\begin{center}
\centerline{\includegraphics[width=1\columnwidth]{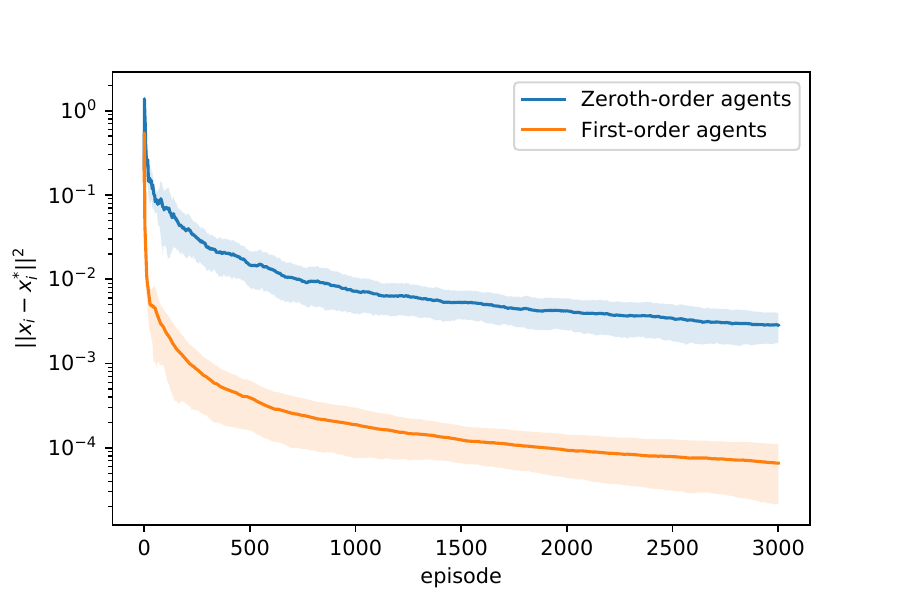}}
\caption{Average error to Nash equilibrium point of first-order agents and zeroth-order agents achieved by Algorithm~\ref{alg:algorithm_asym} when $N_z=5$. The solid lines and shades are averages and standard deviations over 50 runs.}
\label{fig_comp_FOZO_agent}
\end{center}
\end{figure}

\begin{figure}[tb]
\begin{center}
\centerline{\includegraphics[width=1\columnwidth]{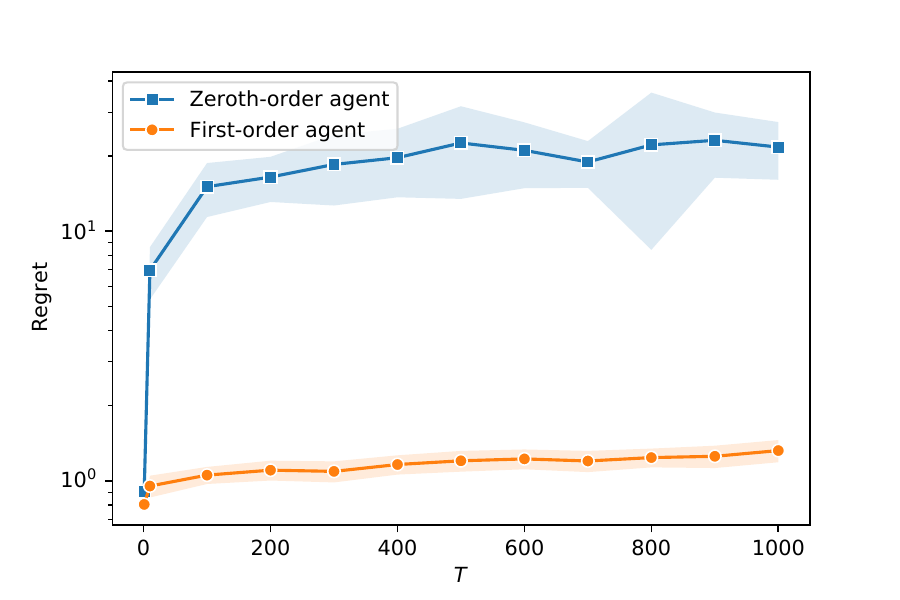}}
\caption{Regret $\sum_{t=1}^T \mathcal{C}_i(x_t) - \mathop{\rm{min}}_{x_i \in \mathcal{X}_i} \sum_{t=1}^T\mathcal{C}_i(x_i,x_{-i,t})$ of first-order agents and zeroth-order agents  achieved by Algorithm~\ref{alg:algorithm_asym} when $N_z=5$. The solid lines and shades are averages and standard deviations over 50 runs.}
\label{fig_regret}
\end{center}
\end{figure}

We run Algorithm~\ref{alg:algorithm_asym} with different values of $N_z$, as well as the pure first- and zeroth-order methods. We choose the step sizes $\eta_{f,t}=\eta_{z,t}=\frac{0.6}{t}$ and the parameter $\delta_t=\frac{0.5}{t}$. Fig.~\ref{fig_NE} illustrates the convergence rate of these algorithms. We observe that our asymmetric learning algorithm with different values of $N_z$ always performs worse than the pure first-order method but better than the pure zeroth-order method, which agrees with Theorem~\ref{theorem:NE}. Besides, a larger value of $N_z$, which means that fewer agents have access to first-order gradient information, leads to a slower convergence rate; when $N_z$ approaches $N$, the convergence rate is similar to that of the pure zeroth-order method.

Fig.~\ref{fig_comp_FOZO_agent} illustrates the error between the action $x_i$ and the corresponding  Nash equilibrium $x^{*}_i$ for every agent $i$. We term the agents with access to first-order gradient information and zeroth-order oracles as first-order agents and zeroth-order agents, respectively. We select $N_z=5$ and compute the average error to the Nash equilibrium for first- and zeroth-order agents. 
As shown in Fig.~\ref{fig_comp_FOZO_agent}, first-order agents converge faster than the zeroth-order agents.
Fig.~\ref{fig_regret} plots the average regret over $50$ runs for first-order agents and zeroth-order agents. We observe that agents that use first-order gradient have much smaller regret than the agents with zeroth-order oracles, which aligns with our theoretical result.

\subsection{Risk-averse Cournot Game}

\begin{figure}[tb]
\begin{center}
\centerline{\includegraphics[width=1\columnwidth]{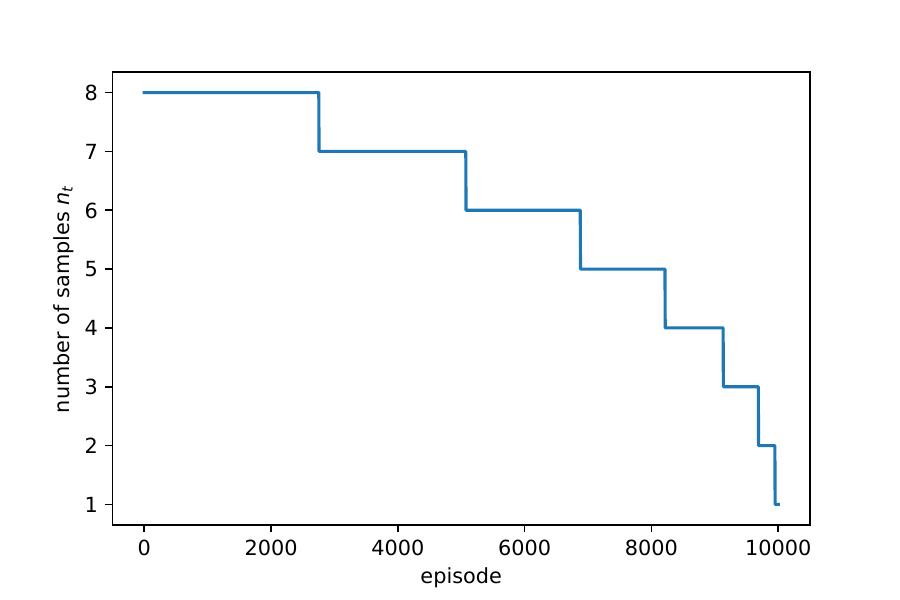}}
\caption{The number of samples of Algorithm~\ref{alg:algorithm_asym} and the pure zeroth-order method (Algorithm 1 in \cite{wang2022risk}).}
\label{fig_number_sample}
\end{center}

\end{figure}

\begin{figure}[tb]
\begin{center}
\centerline{\includegraphics[width=1\columnwidth]{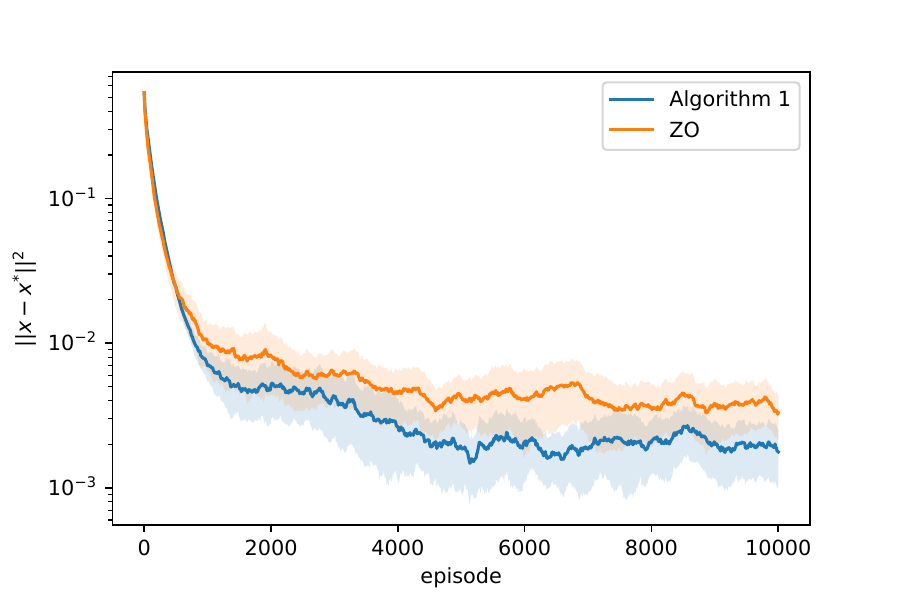}}
\caption{Error to Nash equilibrium point achieved by Algorithm~\ref{alg:algorithm_asym} and the pure zeroth-order method (ZO) in risk-averse games. The solid lines and shades are averages and standard deviations over 50 runs.}
\label{fig_NE_risk}
\end{center}

\end{figure}

Consider a market with two agents, i.e., $N=2$. Each agent decides the quantity $x_i$, $i=1,2$. The stochastic cost of each agent is defined by $J_i(x)= -(2-\sum_{j=1}^2 x_j)x_i+0.2x_i+\xi_i x_i+1$,  where $\xi_i$ is a random variable used to represent the uncertainty in the market. 
Agents aim to minimize the risk of incurring high costs, i.e., agents are risk-averse. We use CVaR as the risk measure and denote the risk level of each agent~$i$ by $\alpha_i$. It is well known that the CVaR value represents the average of the worst $\alpha_i$ percent of the stochastic cost, and when $\alpha_i=1$, it is equivalent to the risk-neutral case.
Let agent~$1$ select the risk level $\alpha_1 = 0.5$ while agent~$2$ is risk-neutral, i.e., $\alpha_2 =1$. The objective functions of these two agents are $\mathcal{C}_1 (x):= {\rm{CVaR}}_{\alpha_1}[J_1(x,\xi_1)]$ and $ \mathcal{C}_2 (x):= \mathbb{E}[J_2(x,\xi_2)]$, respectively, where the definition of CVaR can be found in \cite{wang2022risk}.

The gradient of the CVaR function is hard to compute  in most cases even when the function $J_1(x,\xi_1)$ is completely known \cite{rosen1965existence}. We assume that given the played action profile $x$, the risk-averse agent does not have access to  gradient information but only  to the cost evaluation $J_1(x,\xi_1)$. The risk-neutral agent has access to  gradient information $\nabla_2 \mathcal{C}_2(x)$.

Since the CVaR value cannot be estimated from only one sample, we use the sampling strategy proposed in \cite{wang2022risk} that uses a decreasing number of samples with respect to the number of iterations. The number of samples $n_t$ chosen here is shown in Fig.~\ref{fig_number_sample}. With this sampling strategy, the agents have $n_t$ samples at episode $t$ to estimate the CVaR values. 

Let $\xi_i \sim U(0,0.4)$ be a uniform random variable. Then, we can obtain an explicit expression of the cost function $\mathcal{C}_1(x)$. Further, using the first-order optimality condition we can compute the Nash equilibrium  to be $(0.4667,0.5667)$. We run Algorithm~\ref{alg:algorithm_asym} and the pure zeroth-order algorithm (Algorithm 1 in \cite{wang2022risk}), and select the step size $\eta_t = \frac{0.005}{t}$ and the parameter $\delta = \frac{0.5}{t}$ for both algorithms. Fig.~\ref{fig_NE_risk} shows that our asymmetric learning algorithm converges faster than the pure zeroth-order method.

\section{Conclusion}\label{sec:conclusion}
In this work, we proposed an asymmetric learning algorithm for convex games in which the agents update their actions  using either first-order gradient information or zeroth-order oracles. 
We showed that our algorithm  achieves sub-linear regret for convex games and last-iterate Nash equilibrium convergence for the class of strongly monotone games.
Our theoretical analysis further established that the performance of the proposed algorithm consistently lies between those of the pure first-order and zeroth-order methods. This result highlights the algorithm's flexibility, as it can interpolate between these two extremes by varying the number of agents with access to gradient or oracle information. In particular, our algorithm recovers both the pure first-order and pure zeroth-order methods as special cases.
We demonstrated the effectiveness of our algorithm through numerical experiments on deterministic and stochastic risk-averse Cournot games. These simulations illustrate the adaptability of the algorithm in handling heterogeneous agent information structures.

Several promising directions emerge for future research. One avenue is the development of a generalized theoretical framework to systematically study asymmetric learning algorithms that enable agents to employ multiple strategies simultaneously.
Another intriguing direction involves exploring the fairness in such asymmetric settings.
Besides, extending our analysis to more complex game structures, such as non-convex or hierarchical games, and studying the implications on equilibrium computation and convergence rates would also be interesting.





%



\bibliography{00_citation}
\bibliographystyle{unsrt}







\clearpage

\appendix\label{sec:appendix}
\subsection{Proof of Lemma \ref{lemma:func:property}}\label{appen:fundamental:lemma}

\noindent 1)  Since the function $\mathcal{C}_i(x)$ is convex in $x_i$ for any fixed $x_{-i}$, we have that $\mathcal{C}_{i_z,t}(x_{z})=\mathcal{C}_{i_z}(x_{z},x_{f,t})$
is convex in $x_{i_z}$. Therefore, for any $p_1,p_2 \in \mathcal{X}_{i_z}^{\delta_t}$ and $\theta \in [0,1]$, we have that $\mathcal{C}_{i_z,t}(\theta p_1+ (1-\theta) p_2,x_{-i_z})\leq \theta \mathcal{C}_{i_z,t}( p_1,x_{-i_z}) +(1-\theta) \mathcal{C}_{i_z,t}( p_2,x_{-i_z})$. Based on this, it gives that
\begin{align*}
&\mathcal{C}_{i_z,t}^{\delta_t}(\theta p_1+ (1-\theta) p_2,x_{-i_z})\nonumber \\
&= \mathbb{E} \big[ \mathcal{C}_{i_z,t}(\theta p_1+ (1-\theta) p_2+\delta_t w_{i_z},x_{-i_z}+\delta_t u_{-i_z}) \big] \nonumber \\
&= \mathbb{E} \big[ \mathcal{C}_{i_z,t}(\theta (p_1+\delta_t w_{i_z}) +(1-\theta) (p_2+\delta_t w_{i_z}), \nonumber \\
&\quad x_{-i_z}+\delta_t u_{-i_z}) \big] \nonumber \\
&\leq   \mathbb{E} \big[ \theta  \mathcal{C}_{i_z,t}( p_1+\delta_t w_{i_z},x_{-i_z}+\delta_t u_{-i_z})\nonumber \\
&\quad + (1-\theta)\mathcal{C}_{i_z,t}( p_2+\delta_t w_{i_z},x_{-i_z}+\delta_t u_{-i_z})  \big] \nonumber \\
&=  \theta  \mathcal{C}_{i_z,t}^{\delta_t}( p_1,x_{-i_z}) + (1-\theta)\mathcal{C}_{i_z,t}^{\delta_t}( p_2,x_{-i_z}), \nonumber \\
\end{align*}
where the expectations above are taken with respect to $w_{i_z}\sim \mathbb{B}$ and $u_{-i_z}\sim \mathbb{S}_{-i_z}$. The proof is complete.

\noindent 2) Since the function $\mathcal{C}_i(x)$ is $L_0$-Lipschitz continuous in $x$, we have that $\mathcal{C}_{i_z,t}(x_z)$ is $L_0$-Lipschitz continuous in $x_z$. Thus, for arbitrary two points $x_z, x_z' \in \mathcal{X}_z$ with $\mathcal{X}_z = \Pi_{i_z=1}^{N_z} \mathcal{X}_{i_z}$, we have that 
\begin{align*}
    &\Big|\mathcal{C}_{i_z,t}^{\delta_t}(x_z)-\mathcal{C}_{i_z,t}^{\delta_t}(x_z')\Big|\\
    &=\Big|\mathop{\mathbb{E}}_{w_{i_z}\sim \mathbb{B}, u_{-i_z}\sim \mathbb{S}_{-i_z}}[ \mathcal{C}_{i_z,t}(x_{i_z}+\delta_t w_{i_z},x_{-i_z}+\delta_t u_{-i_z}) \\
    &\quad - \mathcal{C}_{i_z,t}(x_{i_z}'+\delta_t w_{i_z},x_{-i_z}'+\delta_t u_{-i_z}) ]\Big| \\
    &\leq \mathop{\mathbb{E}}_{w_{i_z}\sim \mathbb{B}, u_{-i_z}\sim \mathbb{S}_{-i_z}}[ L_0 \left\| x_z-x_z'\right\|] \\
    & \leq L_0 \left\| x_z-x_z'\right\|.
\end{align*} The proof is complete.

\noindent 3) From the definition of $\mathcal{C}_{i_z,t}$, it directly follows that
\begin{align*}
    &\Big|\mathcal{C}_{i_z,t}^{\delta_t}(x_z)-\mathcal{C}_{i_z,t}(x_z)\Big| \\
    &=\Big| \mathop{\mathbb{E}}_{w_{i_z}\sim \mathbb{B}, u_{-i_z}\sim \mathbb{S}_{-i_z}}[ \mathcal{C}_{i_z,t}(x_{i_z}+\delta_t w_{i_z},x_{-i_z}+\delta_t u_{-i_z}) \\
    &\quad - \mathcal{C}_{i_z,t}(x_z) ] \Big| \\
    &\leq \mathop{\mathbb{E}}_{w_{i_z}\sim \mathbb{B}, u_{-i_z}\sim \mathbb{S}_{-i_z}} L_0  \left\|   (\delta_t w_{i_z},\delta_t u_{-i_z})\right\|  \\
    & \leq L_0 \delta_t \sqrt{N_z}.
\end{align*}

\noindent 4) The proof of this claim can be directly adapted from Lemma~C.1 in \cite{bravo2018bandit} and thus omitted.

\subsection{Proof of Theorem \ref{theorem:No-regret}}\label{appen:theorem:no_regret}
We first focus on the regret analysis of agents with first-order gradient information. 
As a reminder, we use $i_f$ to index an agent with access to first-order gradient information and $i_z$ to index an agent using zeroth-order oracles.
Define $y_{i_f}^{*}:= \mathop{\rm{argmin}}_{x_{i_f} \in \mathcal{X}_{i_f} } \sum_{t=1}^T \mathcal{C}_{i_f}(x_{i_f}, x_{-i_f,t},\hat{x}_{z,t})$. From the update rule \eqref{eq:asy:update:FO}, we have that
\begin{align}\label{eq:no_regret:FO:temp1}
    &\left\| x_{i_f,t+1} - y_{i_f}^{*}\right\|^2 \nonumber \\
    &= \left\| \mathcal{P}_{\mathcal{X}_{i_f}} \Big( x_{i_f,t} - \eta_{f,t} \nabla_{i_f} \mathcal{C}_{i_f}(x_{f,t},\hat{x}_{z,t}) \Big) - y_{i_f}^{*}\right\|^2 \nonumber \\
    &\leq  \left\|  x_{i_f,t} - \eta_{f,t} \nabla_{i_f} \mathcal{C}_{i_f}(x_{f,t},\hat{x}_{z,t}) - y_{i_f}^{*}\right\|^2 \nonumber \\
    &\leq  \left\|  x_{i_f,t}- y_{i_f}^{*}\right\|^2 + \eta_{f,t}^2 \left\| \nabla_{i_f} \mathcal{C}_{i_f}(x_{f,t},\hat{x}_{z,t})\right\|^2 \nonumber \\
    &\quad - 2 \eta_{f,t} \left\langle \nabla_{i_f} \mathcal{C}_{i_f}(x_{f,t},\hat{x}_{z,t}), x_{i_f,t}- y_{i_f}^{*} \right\rangle ,
\end{align}
where the first inequality is due to the facts that $\mathcal{P}_{\mathcal{X}_{i_f}}(y_{i_f}^{*}) = y_{i_f}^{*}$ and the projection operator is non-expansive. 
From the convexity of the function $\mathcal{C}_{i_f}(x)$ in $x_{i_f}$, we have that
\begin{align}
    &{\rm{R}}_{i_f}(T)= \mathbb{E} \Big[ \sum_{t=1}^T \mathcal{C}_{i_f}(x_{f,t},\hat{x}_{z,t}) - \sum_{t=1}^T \mathcal{C}_{i_f}(y_{i_f}^{*}, x_{-i_f,t},\hat{x}_{z,t}) \Big] \nonumber \\
    &\leq  \mathbb{E} \Big[\sum_{t=1}^T \left\langle \nabla_{i_f} \mathcal{C}_{i_f}(x_{f,t},\hat{x}_{z,t}),x_{i_f,t} - y_{i_f}^{*} \right\rangle \Big]\nonumber \\
    &\leq  \sum_{t=1}^T \frac{1}{2\eta_{f,t}}\Big( \mathbb{E} \left\| x_{i_f,t} - y_{i_f}^{*}\right\|^2 - \mathbb{E} \left\| x_{i_f,t+1} - y_{i_f}^{*}\right\|^2  \nonumber \\
    & \quad + \eta_{f,t}^2 \mathbb{E} \left\| \nabla_{i_f} \mathcal{C}_{i_f}(x_{f,t},\hat{x}_{z,t})\right\|^2  \Big) \nonumber \\
    &\leq \sum_{t=1}^T  \mathbb{E} \left\| x_{i_f,t} - y_{i_f}^{*}\right\|^2  \Big(\frac{1}{2\eta_{f,t}} - \frac{1}{2\eta_{f,t-1}} \Big) + \sum_{t=1}^T\frac{\eta_{f,t} L_0^2}{2} \nonumber \\
    &\leq \frac{D^2}{2}\sum_{t=1}^T \Big(\frac{1}{\eta_{f,t}} - \frac{1}{\eta_{f,t-1}} \Big) + \sum_{t=1}^T\frac{\eta_{f,t} L_0^2}{2} \nonumber \\
    &\leq \frac{D^2}{2 \eta_{f,T}}  + \sum_{t=1}^T\frac{\eta_{f,t} L_0^2}{2} \nonumber \\
    &\leq \frac{ \sqrt{T} }{2 } ( D^2 +  2L_0^2)  ,\nonumber
\end{align}
where $\frac{1}{\eta_{f,0}} :=0$. The second inequality follows from \eqref{eq:no_regret:FO:temp1} and the last inequality follows since $\eta_{f,t}=\frac{1}{ \sqrt{t}} $ and $\sum_{t=1}^T \frac{1}{\sqrt{t}} \leq 2 \sqrt{T}$. 
The third inequality holds since, for any two sequences $\{ a_t\}_{t=1}^T$, $\{b_t \}_{t=1}^T$, with $\frac{1}{a_0} :=0$, $a_t>0, b_t\geq0, \forall t=1,\ldots,T$, we have
\begin{align*}
    \sum_{t=1}^T \frac{1}{a_t}(b_t - b_{t+1}) &= \sum_{t=1}^T b_t (\frac{1}{a_t} - \frac{1}{a_{t-1}}) - \frac{b_{T+1}}{a_T} \\&\leq \sum_{t=1}^T b_t (\frac{1}{a_t} - \frac{1}{a_{t-1}}).
\end{align*}
The proof for the agents with first-order gradient information is complete. 

Next we turn to the regret analysis of the agents with zeroth-order oracles. Define  $y_{i_z}^{*}:= \mathop{\rm{argmin}}_{x_{i_z} \in \mathcal{X}_{i_z} } \sum_{t=1}^T \mathcal{C}_{i_z}(x_{i_z}, \hat{x}_{-i_z,t},x_{f,t})$ and $y_{i_z}^{\delta_t}:= (1-\frac{\delta_t}{r})y_{i_z}^{*} $. Since $y_{i_z}^{*} \in \mathcal{X}_{i_z} $, we have $y_{i_z}^{\delta_t} \in \mathcal{X}_{i_z}^{\delta_t}$.

By virtue of the update rule \eqref{eq:asy:update:ZO}, we have that 
\begin{align}\label{eq:no_regret:ZO:temp2}
    &\left\|x_{i_z,t+1} - y_{i_z}^{\delta_t }\right\|^2 \nonumber \\
    &= \left\| \mathcal{P}_{\mathcal{X}_{i_z}^{\delta_t}} \Big( x_{i_z,t} - \eta_{z,t} {g}_{i_z,t}\Big) - y_{i_z}^{\delta_t}  \right\| ^2 \nonumber \\
    &\leq \left\|   x_{i_z,t} - \eta_{z,t} {g}_{i_z,t} - y_{i_z}^{\delta_t}  \right\| ^2\nonumber \\
    &\leq  \left\|x_{i_z,t} - y_{i_z}^{\delta_t}\right\|^2 + \eta_{z,t}^2 \left\| {g}_{i_z,t}\right\|^2 \nonumber \\
    &\quad - 2 \eta_{z,t} \left\langle x_{i_z,t} -y_{i_z}^{\delta_t}, g_{i_z,t}\right\rangle,
\end{align}
where the first inequality is due to the fact that $\mathcal{P}_{\mathcal{X}_{i_z}^{\delta_t}}( y_{i_z}^{\delta_t}) = y_{i_z}^{\delta_t}$. 
Taking expectations of both sides on the inequality \eqref{eq:no_regret:ZO:temp2}, we obtain that
\begin{align}\label{eq:no_regret:ZO:temp3}
    \mathbb{E}\left\|x_{i_z,t+1} - y_{i_z}^{\delta_t}\right\|^2
    \leq  \mathbb{E}\left\|x_{i_z,t} - y_{i_z}^{\delta_t}\right\|^2 + \frac{\eta_{z,t}^2 d^2 U^2}{\delta_t^2}  \nonumber \\
    - 2 \eta_{z,t} \mathbb{E} \left\langle x_{i_z,t} -y_{i_z}^{\delta_t}, \nabla_{i_z} \mathcal{C}_{i_z,t}^{\delta_t}(x_{z,t})\right\rangle,
\end{align}
where the inequality follows from Assumption~\ref{assump:U} and the fact that $\left\| g_{i_z,t}\right\| = \left\| \frac{d}{\delta_t} \mathcal{C}_{i_z}(\hat{x}_{z,t},x_{f,t}) u_{i_z,t}\right\| \leq \frac{d U}{\delta_t}$.
Then, we have that
\begin{align}\label{eq:no_regret:ZO:temp4}
    &{\rm{R}}_{i_z}(T) = \mathbb{E} \Big[\sum_{t=1}^T \mathcal{C}_{i_z} (\hat{x}_{z,t},x_{f,t}) - \sum_{t=1}^T \mathcal{C}_{i_z} (y_{i_z}^{*},\hat{x}_{-i_z,t},x_{f,t}) \Big] \nonumber \\
    & = \mathbb{E} \Big[ \sum_{t=1}^T \Big( \mathcal{C}_{i_z,t} (\hat{x}_{z,t}) - \mathcal{C}_{i_z,t} (y_{i_z}^{*},\hat{x}_{-i_z,t}) - \mathcal{C}_{i_z,t} ({x}_{z,t})  \nonumber \\ 
    &\quad +  \mathcal{C}_{i_z,t} ({x}_{z,t}) -\mathcal{C}_{i_z,t} (y_{i_z}^{*},{x}_{-i_z,t}) + \mathcal{C}_{i_z,t} (y_{i_z}^{*},{x}_{-i_z,t})\Big) \Big]\nonumber \\
    &\leq \mathbb{E} \Big[\sum_{t=1}^T \Big( \mathcal{C}_{i_z,t} ({x}_{z,t}) - \mathcal{C}_{i_z,t} (y_{i_z}^{*},{x}_{-i_z,t}) +2 L_0 \delta_t \sqrt{N_z} \Big) \Big] \nonumber \\
    &\leq \mathbb{E} \Big[ \sum_{t=1}^T \Big( \mathcal{C}_{i_z,t}^{\delta_t} ({x}_{z,t}) - \mathcal{C}_{i_z,t}^{\delta_t} (y_{i_z}^{*},{x}_{-i_z,t}) +4 L_0 \delta_t \sqrt{N_z}    \Big) \Big] \nonumber \\
    & = \mathbb{E} \Big[\sum_{t=1}^T \Big( \mathcal{C}_{i_z,t}^{\delta_t} ({x}_{z,t}) - \mathcal{C}_{i_z,t}^{\delta_t} (y_{i_z}^{\delta_t},{x}_{-i_z,t}) + \mathcal{C}_{i_z,t}^{\delta_t} (y_{i_z}^{\delta_t},{x}_{-i_z,t}) \nonumber \\
    &\quad - \mathcal{C}_{i_z,t}^{\delta_t} (y_{i_z}^{*},{x}_{-i_z,t}) +4 L_0 \delta_t \sqrt{N_z}    \Big) \Big] \nonumber \\
    &\leq \mathbb{E} \Big[\sum_{t=1}^T \Big( \mathcal{C}_{i_z,t}^{\delta_t} ({x}_{z,t}) - \mathcal{C}_{i_z,t}^{\delta_t} (y_{i_z}^{\delta_t},{x}_{-i_z,t}) + \frac{L_0 \delta_t  D}{r} \nonumber \\
    &\quad +4 L_0 \delta_t \sqrt{N_z}    \Big) \Big] \nonumber \\
    &\leq  \sum_{t=1}^T \mathbb{E} \left\langle x_{i_z,t} -y_{i_z}^{\delta_t}, \nabla_{i_z} \mathcal{C}_{i_z,t}^{\delta}(x_{z,t})\right\rangle   \nonumber \\
    & \quad + L_0( \frac{D}{r} + 4 \sqrt{N_z}) \sum_{t=1}^T \delta_t.
\end{align}
The first inequality follows from  Assumption~\ref{assump:L0} and the second inequality follows from the third item in Lemma \ref{lemma:func:property}. The third inequality holds because 
$$|\mathcal{C}_{i_z,t}^{\delta_t} (y_{i_z}^{\delta_t},{x}_{-i_z,t}) - \mathcal{C}_{i_z,t}^{\delta_t} (y_{i_z}^{*},{x}_{-i_z,t})| \leq L_0 \frac{\delta_t}{r} \left\| y_{i_z}^{*}\right\| \leq \frac{L_0 \delta_t D}{r}.$$
The last inequality in \eqref{eq:no_regret:ZO:temp4} follows from the first item in Lemma \ref{lemma:func:property}. 
Define $\frac{1}{\eta_{z,0}}=0$.
Substituting the inequality \eqref{eq:no_regret:ZO:temp3}  
into the inequality \eqref{eq:no_regret:ZO:temp4}, we have that
\begin{align}\label{eq:no_regret:ZO:temp5}
    &{\rm{R}}_{i_z}(T)  \nonumber \\
    &\leq  \sum_{t=1}^T \frac{1}{2 \eta_{z,t}} \Big( \mathbb{E}\left\|x_{i_z,t} - y_{i_z}^{\delta_t}\right\|^2 - \mathbb{E}\left\|x_{i_z,t+1} - y_{i_z}^{\delta_t}\right\|^2\Big) \nonumber \\
    &\quad + \sum_{t=1}^T \frac{\eta_{z,t} d^2 U^2}{2 \delta_t^2} + L_0( \frac{D}{r} + 4 \sqrt{N_z}) \sum_{t=1}^T \delta_t  \nonumber \\
    & = \sum_{t=1}^T \frac{1}{2 \eta_{z,t}} \Big( \mathbb{E}\left\|x_{i_z,t} - y_{i_z}^{\delta_t}\right\|^2 - \mathbb{E}\left\|x_{i_z,t+1} - y_{i_z}^{\delta_{t+1}}\right\|^2 \nonumber \\
    &\quad + \mathbb{E}\left\|x_{i_z,t+1} - y_{i_z}^{\delta_{t+1}}\right\|^2 - \mathbb{E}\left\|x_{i_z,t+1} - y_{i_z}^{\delta_t}\right\|^2\Big) \nonumber \\
    &\quad + \sum_{t=1}^T \frac{\eta_{z,t} d^2 U^2}{2 \delta_t^2}+ L_0( \frac{D}{r} + 4 \sqrt{N_z}) \sum_{t=1}^T \delta_t  \nonumber \\
    &\leq \frac{1}{2}\sum_{t=1}^T \mathbb{E}\left\|x_{i_z,t} - y_{i_z}^{\delta_t}\right\|^2 (\frac{1}{ \eta_{z,t}} - \frac{1}{ \eta_{z,t-1}}) \nonumber \\
    & \quad + \sum_{t=1}^T  \frac{1}{2 \eta_{z,t}} \mathbb{E} \left[ \langle y_{i_z}^{\delta_{t+1}} - y_{i_z}^{\delta_{t}} , y_{i_z}^{\delta_t} +y_{i_z}^{\delta_{t+1}} - 2x_{i_z,t+1} \rangle  \right]\nonumber \\
    &\quad + \sum_{t=1}^T \frac{\eta_{z,t} d^2 U^2}{2 \delta_t^2}+ L_0( \frac{D}{r} + 4 \sqrt{N_z}) \sum_{t=1}^T \delta_t  \nonumber \\
    &\leq  \frac{D^2}{2 \eta_{z,T}} + \sum_{t=1}^T \frac{\eta_{z,t} d^2 U^2}{2 \delta_t^2}+  L_0( \frac{D}{r} + 4 \sqrt{N_z}) \sum_{t=1}^T \delta_t \nonumber \\
    &\quad + \sum_{t=1}^T \frac{1}{2 \eta_{z,t}} \mathbb{E} \left\|  y_{i_z}^{\delta_{t+1}} - y_{i_z}^{\delta_{t}}\right\|  2D \nonumber \\
    &\leq  \frac{D^2}{2 \eta_{z,T}} + \sum_{t=1}^T \frac{\eta_{z,t} d^2 U^2}{2 \delta_t^2} +  L_0( \frac{D}{r} + 4 \sqrt{N_z}) \sum_{t=1}^T \delta_t \nonumber \\
    &\quad +  \sum_{t=1}^T \frac{(\delta_t-\delta_{t+1})D^2}{r \eta_{z,t}} \nonumber \\
    &\leq  \frac{D^2}{2 \eta_{z,T}} + \sum_{t=1}^T \frac{\eta_{z,t} d^2 U^2}{2 \delta_t^2} +  L_0( \frac{D}{r} + 4 \sqrt{N_z}) \sum_{t=1}^T \delta_t \nonumber \\ 
    &\quad +  \frac{D^2}{r} \sum_{t=1}^T  \delta_t \Big( \frac{1}{\eta_{z,t}}  - \frac{1}{\eta_{z,t-1}} \Big).
\end{align}
In the second and last inequalities, we use again the trick that $\sum_{t=1}^T \frac{1}{a_t}(b_t - b_{t+1}) \leq \sum_{t=1}^T b_t (\frac{1}{a_t} - \frac{1}{a_{t-1}})$ for any two sequences $\{ a_t\}_{t=1}^T$ and $\{ b_t\}_{t=1}^T$ with $\frac{1}{a_0} :=0$, $a_t>0, b_t\geq0, \forall t=[1,T]$. 
The fourth inequality holds since $\left\|y_{i_z}^{\delta_{t+1}} - y_{i_z}^{\delta_t}\right\| =\left\| \frac{(\delta_{t+1}-\delta_t )y_{i_z}^{*}}{r} \right\| \leq \frac{(\delta_t -\delta_{t+1} ) D}{r}$.
Substituting $\eta_{z,t} = \frac{1}{N_z^{\frac{1}{4}} t^{\frac{3}{4}}}$ and $\delta_t =  \frac{1}{N_z^{\frac{1}{4}} t^{\frac{1}{4}}} $ into the inequality \eqref{eq:no_regret:ZO:temp5}, we have
\begin{align*}
    &{\rm{R}}_{i_z}(T)  \nonumber \\
    &\leq \frac{D^2 N_z^{\frac{1}{4}} T^{\frac{3}{4}}}{2} +  \frac{L_0( \frac{D}{r} + 4 \sqrt{N_z})}{N_z^{\frac{1}{4}} } \sum_{t=1}^T \frac{1}{ t^{\frac{1}{4}}} +\frac{d^2 U^2}{2}\sum_{t=1}^T\frac{1}{t}  \nonumber \\
    &\quad + \frac{D^2}{r} \sum_{t=1}^T \frac{1}{ t^{\frac{1}{4}}} \big( t^{\frac{3}{4}} - (t-1)^{\frac{3}{4}}\big) \nonumber \\
    &\leq \frac{D^2 N_z^{\frac{1}{4}} T^{\frac{3}{4}}}{2} +  \frac{L_0( \frac{D}{r} + 4 \sqrt{N_z})}{N_z^{\frac{1}{4}} } \sum_{t=1}^T \frac{1}{ t^{\frac{1}{4}}} + \frac{d^2 U^2 \ln T}{2} \nonumber \\
    &\quad + \frac{D^2}{r} \sum_{t=1}^T \frac{1}{ t^{\frac{1}{4}}} \nonumber \\
    &\leq \frac{D^2 N_z^{\frac{1}{4}} T^{\frac{3}{4}}}{2} +  \frac{4L_0( \frac{D}{r} + 4 \sqrt{N_z})}{3N_z^{\frac{1}{4}} } T^{\frac{3}{4}} + \frac{d^2 U^2 \ln T}{2} \nonumber \\
    &\quad + \frac{4D^2}{3} T^{\frac{3}{4}}\nonumber \\
    &=  \mathcal{O}\left( N_Z^{\frac{1}{4}} T^{\frac{3}{4}} \right),
\end{align*}
where the last inequality follows from the fact that $\sum_{t=1}^T  t^{-\frac{1}{4}} \leq 1+ \int_{1}^T t^{-\frac{1}{4}} dt \leq \frac{4}{3} T^{\frac{3}{4}}$.  The second inequality holds since $ f(t):= t^{\frac{3}{4}} - (t-1)^{\frac{3}{4}} \leq 1$ for all $t\geq 1$, which is easily verified through computing the derivative of $f(t)$.
The proof is complete.

\subsection{Proof of Theorem \ref{theorem:NE}}\label{appen:theorem:NE}

To clarify the notations, we define the following:
\begin{align*}
    x^{*}& = (x_z^{*},x_f^{*}),  \; x_z^{*}= (x_{1}^{*},x_2^{*},\ldots, x_{N_z}^{*}), \\
    x_f^{*} &= (x_{N_z+1}^{*},x_{N_z+2}^{*}, \ldots, x_N^{*}), \; x_t = (x_{z,t},x_{f,t}),\\
    x_{z,t} & = (x_{1,t},\ldots, x_{N_z,t}), \; x_{f,t}= (x_{N_z+1,t},\ldots, x_{N,t}) ,\\
    x_t^{*}& = (x_{z,t}^{*}, x_f^{*}), \; x_{z,t}^{*} = (1-\frac{\delta_t}{r})x_z^{*}, \\
    A_{i_z,t} :&= \left\| x_{i_z,t} - x_{i_z}^{*}\right\|^2 , \; A_{i_f,t} := \left\| x_{i_f,t} - x_{i_f}^{*}\right\|^2,\\
    \bar{A}_{i_z,t} :&= \left\| x_{i_z,t} - x_{i_z,t-1}^{*}\right\|^2 , \; \hat{A}_{i_z,t} := \left\| x_{i_z,t} - x_{i_z,t}^{*}\right\|^2,\\
    A_t :&= \sum_{i_z=1}^{N_z} A_{i_z,t} + \sum_{i_f=N_z+1}^{N} A_{i_f,t} = \left\|x_t - x^{*} \right\|^2, \nonumber \\
    \bar{A}_t :&= \sum_{i_z=1}^{N_z} \bar{A}_{i_z,t} + \sum_{i_f=N_z+1}^{N} A_{i_f,t} = \left\|x_t - x_{t-1}^{*} \right\|^2, \nonumber \\
    \hat{A}_t :&= \sum_{i_z=1}^{N_z} \hat{A}_{i_z,t} + \sum_{i_f=N_z+1}^{N} A_{i_f,t} = \left\|x_t - x_{t}^{*} \right\|^2.
\end{align*}
Recall that we use $i_f$ to index an agent with access to first-order gradient information and $i_z$ to index an agent using zeroth-order oracles.
We first focus on the evolution of $A_{i_f,t}$, i.e., the error dynamics of the agents with first-order gradient information. From the update rule \eqref{eq:asy:update:FO}, we have that
\begin{align}\label{eq:sec3:temp1}
    &A_{i_f,t+1} =\left\| x_{i_f,t+1} - x_{i_f}^{*}\right\|^2 \nonumber \\
    &= \left\| \mathcal{P}_{\mathcal{X}_{i_f}} \Big( x_{i_f,t} - \eta_{f,t} \nabla_{i_f} \mathcal{C}_{i_f}(x_{f,t},\hat{x}_{z,t}) \Big) - x_{i_f}^{*}\right\|^2 \nonumber \\
    &\leq  \left\|  x_{i_f,t} - \eta_{f,t} \nabla_{i_f} \mathcal{C}_{i_f}(x_{f,t},\hat{x}_{z,t}) - x_{i_f}^{*}\right\|^2 \nonumber \\
    &\leq \left\|  x_{i_f,t}- x_{i_f}^{*}\right\|^2 + \eta_{f,t}^2 \left\| \nabla_{i_f} \mathcal{C}_{i_f}(x_{f,t},\hat{x}_{z,t})\right\|^2 \nonumber \\
    &\quad - 2 \eta_{f,t} \left\langle \nabla_{i_f} \mathcal{C}_{i_f}(x_{f,t},\hat{x}_{z,t}), x_{i_f,t}- x_{i_f}^{*} \right\rangle \nonumber \\
    &\leq  A_{i_f,t} - 2 \eta_{f,t} \left\langle \nabla_{i_f} \mathcal{C}_{i_f}(x_{f,t},\hat{x}_{z,t}), x_{i_f,t}- x_{i_f}^{*} \right\rangle \nonumber \\
    &\quad + \eta_{f,t}^2 L_0^2,
\end{align}
where the first inequality follows since $\mathcal{P}_{\mathcal{X}_{i_f}}(x_{i_f}^{*}) = x_{i_f}^{*}$ and the projection operator is non-expansive.
Since $x^{*}$ is a Nash equilibrium of the convex game \eqref{eq:def:game}, we have that $ \langle  \nabla_{i_f} \mathcal{C}_{i_f}(x^{*}), x_{i_f,t}-x_{i_f}^{*}\rangle \geq 0$, $i=1,\ldots,N$. Combining this Nash equilibrium condition with the inequality \eqref{eq:sec3:temp1}, we have that
\begin{align}\label{eq:sec3:temp2}
    &A_{i_f,t+1} \nonumber \\
    &\leq A_{i_f,t} + \eta_{f,t}^2 L_0^2 +  2 \eta_{f,t}  \left\langle \nabla_{i_f} \mathcal{C}_{i_f}(x^{*}),x_{i_f,t}- x_{i_f}^{*} \right\rangle \nonumber \\
    & \quad - 2 \eta_{f,t} \left\langle \nabla_{i_f} \mathcal{C}_{i_f}(x_{f,t} ,\hat{x}_{z,t}) , x_{i_f,t}- x_{i_f}^{*} \right\rangle \nonumber \\
    &=  A_{i_f,t} - 2 \eta_{f,t} \Big\langle \nabla_{i_f} \mathcal{C}_{i_f}(x_{f,t},\hat{x}_{z,t}) -\nabla_{i_f} \mathcal{C}_{i_f}(x_t),\nonumber \\
    &\quad \; x_{i_f,t}- x_{i_f}^{*} \Big\rangle + \eta_{f,t}^2 L_0^2  \nonumber \\
    & \quad -  2 \eta_{f,t} \left\langle \nabla_{i_f} \mathcal{C}_{i_f}(x_t)- \nabla_{i_f} \mathcal{C}_{i_f}(x_t^{*}),x_{i_f,t}- x_{i_f}^{*} \right\rangle \nonumber \\
    & \quad -  2 \eta_{f,t} \left\langle \nabla_{i_f} \mathcal{C}_{i_f}(x_t^{*})- \nabla_{i_f} \mathcal{C}_{i_f}(x^{*}),x_{i_f,t}- x_{i_f}^{*} \right\rangle \nonumber \\
    &\leq  A_{i_f,t} + \eta_{f,t}^2 L_0^2 + 2 \eta_{f,t} L_1 \delta_t \sqrt{N_z} D + \frac{2 \eta_{f,t} L_1 \delta_t D ^2 \sqrt{N_z}}{r}\nonumber \\
    & \quad -  2 \eta_{f,t} \left\langle \nabla_{i_f} \mathcal{C}_{i_f}(x_t)- \nabla_{i_f} \mathcal{C}_{i_f}(x_t^{*}),x_{i_f,t}- x_{i_f}^{*} \right\rangle ,
\end{align}
where the last inequality follows from Assumptions~\ref{assump:L0}, \ref{assump:L1}, and the facts that $\left\| \nabla_{i_f} \mathcal{C}_{i_f}(x_{f,t},\hat{x}_{z,t}) -\nabla_{i_f} \mathcal{C}_{i_f}(x_t)\right\| \leq L_1 \left\| \hat{x}_{z,t} - x_{z,t}\right\| \leq L_1 \delta_t \sqrt{N_z}$, $ \left\| \nabla_{i_f} \mathcal{C}_{i_f}(x_t^{*})- \nabla_{i_f} \mathcal{C}_{i_f}(x^{*})\right\| \leq L_1 \left\| x_t^{*} - x^{*} \right\| = L_1 \left\| x_{z,t}^{*} - x_z^{*} \right\| \leq \frac{L_1 \delta_t D \sqrt{N_z}}{r}$. 

Next we turn our attention on the evolution of $\bar{A}_{i_z,t}$. Similarly, from the update equation \eqref{eq:asy:update:ZO}, we obtain that
\begin{align}\label{eq:sec3:temp3}
    &\bar{A}_{i_z,t+1} =\left\| x_{i_z,t+1} - x_{i_z,t}^{*}\right\|^2 \nonumber \\
    &= \left\| \mathcal{P}_{\mathcal{X}_{i_z}^{\delta_t}} \Big( x_{i_z,t} - \eta_{z,t} g_{i_z,t} \Big) - x_{i_z,t}^{*}\right\|^2 \nonumber \\
    &\leq  \left\|  x_{i_z,t} - \eta_{z,t} g_{i_z,t} - x_{i_z,t}^{*}\right\|^2 \nonumber \\
    &=  \hat{A}_{i_z,t} + \eta_{z,t}^2 \left\| g_{i_z,t}\right\|^2 -   2 \eta_{z,t} \left\langle g_{i_z,t}, x_{i_z,t}- x_{i_z,t}^{*} \right\rangle\nonumber \\
    &\leq  \hat{A}_{i_z,t} + \frac{\eta_{z,t}^2 d^2 U^2}{\delta_t^2}  -  2 \eta_{z,t} \left\langle g_{i_z,t}, x_{i_z,t}- x_{i_z,t}^{*} \right\rangle,
\end{align}
where the first inequality follows from the fact that $\mathcal{P}_{\mathcal{X}_{i_z}^{\delta_t}} \Big(x_{i_z,t}^{*}\Big) = x_{i_z,t}^{*}$ and the projection operator is non-expansive.
Taking expectation of both sides of the above inequality \eqref{eq:sec3:temp3}, we have that

\begin{align}\label{eq:sec3:temp4}
    &\mathbb{E}[ \bar{A}_{i_z,t+1}] \nonumber \\
    &\leq   \mathbb{E}  [\hat{A}_{i_z,t}] -  2 \eta_{z,t} \mathbb{E} \left\langle \nabla_{i_z} \mathcal{C}^{\delta_t}_{i_z,t}(x_{z,t}), x_{i_z,t}- x_{i_z,t}^{*} \right\rangle  \nonumber \\
    &\quad  + \frac{\eta_{z,t}^2 d^2 U^2}{\delta_t^2}  \nonumber \\
    &\leq  \mathbb{E}  [\hat{A}_{i_z,t}] + \frac{\eta_{z,t}^2 d^2 U^2}{\delta_t^2} \nonumber \\
    &\quad -  2 \eta_{z,t} \mathbb{E}\left\langle \nabla_{i_z} \mathcal{C}^{\delta_t}_{i_z,t}(x_{z,t})  - \nabla_{i_z} \mathcal{C}_{i_z}(x_t) , x_{i_z,t}- x_{i_z,t}^{*} \right\rangle \nonumber \\
    & \quad -  2 \eta_{z,t} \mathbb{E}\left\langle \nabla_{i_z} \mathcal{C}_{i_z}(x_t)-\nabla_{i_z} \mathcal{C}_{i_z}(x_t^{*}),x_{i_z,t}- x_{i_z,t}^{*} \right\rangle \nonumber \\
    & \quad -  2 \eta_{z,t} \mathbb{E}\left\langle \nabla_{i_z} \mathcal{C}_{i_z}(x_t^{*})-\nabla_{i_z} \mathcal{C}_{i_z}(x^{*}),x_{i_z,t}- x_{i_z,t}^{*} \right\rangle \nonumber \\
    & \quad -2 \eta_{z,t} \mathbb{E} \left\langle \nabla_{i_z} \mathcal{C}_{i_z}(x^{*}),(x_{i_z,t}- x_{i_z,t}^{*}) - (x_{i_z,t}- x_{i_z}^{*}) \right\rangle \nonumber \\
    &\leq  \mathbb{E}  [\hat{A}_{i_z,t}] + \frac{\eta_{z,t}^2 d^2 U^2}{\delta_t^2} + \frac{2\eta_{z,t} L_1 \delta_t D^2}{r} \nonumber \\
    &\quad -  2 \eta_{z,t} \mathbb{E}\left\langle \nabla_{i_z} \mathcal{C}^{\delta_t}_{i_z,t}(x_{z,t})  - \nabla_{i_z} \mathcal{C}_{i_z}(x_t) , x_{i_z,t}- x_{i_z,t}^{*} \right\rangle \nonumber \\
    & \quad -  2 \eta_{z,t} \mathbb{E}\left\langle \nabla_{i_z} \mathcal{C}_{i_z}(x_t)-\nabla_{i_z} \mathcal{C}_{i_z}(x_t^{*}),x_{i_z,t}- x_{i_z,t}^{*} \right\rangle \nonumber \\
    & \quad  -2 \eta_{z,t} \mathbb{E} \left\langle \nabla_{i_z} \mathcal{C}_{i_z}(x^{*}),x_{i_z}^{*}   - x_{i_z,t}^{*} \right\rangle \nonumber \\
    &\leq  \mathbb{E}  [\hat{A}_{i_z,t}] + \frac{\eta_{z,t}^2 d^2 U^2}{\delta_t^2} + \frac{2\eta_{z,t} L_1 \delta_t D^2}{r} + \frac{2\eta_{z,t} L_0 \delta_t D }{r} \nonumber \\
    &\quad -  2 \eta_{z,t} \mathbb{E}\left\langle \nabla_{i_z} \mathcal{C}^{\delta_t}_{i_z,t}(x_{z,t})  - \nabla_{i_z} \mathcal{C}_{i_z}(x_t) , x_{i_z,t}- x_{i_z,t}^{*} \right\rangle \nonumber \\
    & \quad -  2 \eta_{z,t} \mathbb{E}\left\langle \nabla_{i_z} \mathcal{C}_{i_z}(x_t)-\nabla_{i_z} \mathcal{C}_{i_z}(x_t^{*}),x_{i_z,t}- x_{i_z,t}^{*} \right\rangle, 
\end{align} where 
the first inequality follows from Lemma \ref{lemma:func:property}, the second inequality holds due to the Nash equilibrium first-order optimality condition, and the third inequality follows since $\left\|\nabla_{i_z} \mathcal{C}_{i_z}(x_t^{*})-\nabla_{i_z} \mathcal{C}_{i_z}(x^{*})\right\| \leq L_1 \left\| x_t^{*} - x^{*}\right\| \leq \frac{L_1 \delta_t D}{r}$. The last inequality follows from Assumption~\ref{assump:L0} and the fact that $ \left\| x_{i_z}^{*}   - x_{i_z,t}^{*}\right\| \leq \frac{\delta_t D}{r} $.

Now we are in a position to bound the term $\left\| \nabla_{i_z} C^{\delta_t}_{i_z,t}(x_{z,t})  - \nabla_{i_z} C_{i_z}(x_t)\right\|$. Since $ \mathcal{C}_{i_z}$ is bounded and continuous for all $i_z \in \mathcal{N}_z$, by Lebesgue’s dominated convergence theorem (Chapter 4 in \cite{royden1988real}), the order of integration and differentiation can be interchanged. Then, it follows that
\begin{align*}
    & \nabla_{i_z} C^{\delta_t}_{i_z,t}(x_{z,t})   \nonumber \\
    &=   \nabla_{i_z} \mathbb{E} \left[ C_{i_z,t}(x_{i_z,t} +\delta_t w_{i_z,t}, x_{-i_z,t}+\delta_t u_{-i_z,t})\right]  \nonumber \\
    &=    \mathbb{E} \left[  \nabla_{i_z} C_{i_z,t}(x_{i_z,t} +\delta_t w_{i_z,t}, x_{-i_z,t}+\delta_t u_{-i_z,t}) \right]   \nonumber \\
    &=    \mathbb{E} \left[  \nabla_{i_z} C_{i_z}(x_{i_z,t} +\delta_t w_{i_z,t}, x_{-i_z,t}+\delta_t u_{-i_z,t},x_{f,t}) \right]   \nonumber,
\end{align*} 
where the expectations in the above inequality are taken w.r.t $w_{i_z,t} \sim \mathbb{B}$ and $u_{-i_z,t} \sim \mathbb{S}_{-i_z}$. Then, we have that
\begin{align}\label{eq:sec3:temp5}
    &\left\| \nabla_{i_z} C^{\delta_t}_{i_z,t}(x_{z,t})  - \nabla_{i_z} C_{i_z}(x_t)\right\| \nonumber \\
    &\leq  L_1 \delta_t \left\| ( w_{i_z,t}, u_{-i_z,t})\right\| \leq L_1 \delta_t \sqrt{N_z}.
\end{align} 
Substituting \eqref{eq:sec3:temp5} into \eqref{eq:sec3:temp4}, we obtain that 
\begin{align}\label{eq:sec3:temp6}
    &\mathbb{E}[ \bar{A}_{i_z,t+1}] \nonumber \\
    &\leq  \mathbb{E}  [\hat{A}_{i_z,t}] + \frac{\eta_{z,t}^2 d^2 U^2}{\delta_t^2}  + \frac{2\eta_{z,t} L_1 \delta_t D^2}{r} + \frac{2\eta_{z,t} L_0 \delta_t D }{r}\nonumber \\
    &\quad -  2 \eta_{z,t} \mathbb{E} \left\langle \nabla_{i_z} C_{i_z}(x_t)-\nabla_{i_z} \mathcal{C}_{i_z}(x_t^{*}),x_{i_z,t}- x_{i_z,t}^{*} \right\rangle \nonumber \\
    &\quad + 2 \eta_{z,t} L_1 \delta_t \sqrt{N_z} D.
\end{align}
Define $\delta_0 =0$ and $x_{i_z,0}^{*} = x_{i_z}^{*}$. From the definition of $\hat{A}_{i_z,t}$ and $\bar{A}_{i_z,t}$, we have
\begin{align}\label{eq:sec3:temp6.1}
    &\hat{A}_{i_z,t} = \left\| x_{i_z,t} - x_{i_z,t-1}^{*} +x_{i_z,t-1}^{*}  - x_{i_z,t}^{*}\right\|^2 \nonumber \\
    &= \bar{A}_{i_z,t} + \left\| x_{i_z,t-1}^{*}  - x_{i_z,t}^{*}\right\|^2 \nonumber \\
    &\quad + 2 \langle x_{i_z,t} - x_{i_z,t-1}^{*}, x_{i_z,t-1}^{*}  - x_{i_z,t}^{*} \rangle \nonumber \\
    &\leq \bar{A}_{i_z,t} + \frac{|\delta_{t-1} - \delta_t|^2 D^2}{r^2} + \frac{2 |\delta_{t-1} - \delta_t| D^2}{r},
\end{align}
for all $t\geq1$. The inequality \eqref{eq:sec3:temp6.1} holds at $t=1$ since $ \left\| x_{i_z,0}^{*}  - x_{i_z,1}^{*}\right\| = \frac{\delta_1}{r} x_{i_z}^{*}  \leq \frac{|\delta_0 -\delta_1|D}{r}$.
Substituting \eqref{eq:sec3:temp6.1} into \eqref{eq:sec3:temp6}, we have
\begin{align}\label{eq:sec3:temp6.2}
    &\mathbb{E}[ \bar{A}_{i_z,t+1}] \nonumber \\
    &\leq  \mathbb{E}  [\bar{A}_{i_z,t}] + \frac{|\delta_{t-1} - \delta_t|^2 D^2}{r^2} + \frac{2 |\delta_{t-1} - \delta_t| D^2}{r} \nonumber \\
    &\quad -  2 \eta_{z,t} \mathbb{E} \left\langle \nabla_{i_z} C_{i_z}(x_t)-\nabla_{i_z} \mathcal{C}_{i_z}(x_t^{*}),x_{i_z,t}- x_{i_z,t}^{*} \right\rangle \nonumber \\
    &\quad + \frac{\eta_{z,t}^2 d^2 U^2}{\delta_t^2}  + \frac{2\eta_{z,t} L_1 \delta_t D^2}{r} + \frac{2\eta_{z,t} L_0 \delta_t D }{r}\nonumber \\
    &\quad + 2 \eta_{z,t} L_1 \delta_t \sqrt{N_z} D.
\end{align}
Summing \eqref{eq:sec3:temp2} over $i_f=N_z+1,\ldots,N$ and \eqref{eq:sec3:temp6.2} over $i_z=1,\ldots,N_z$ altogether, and setting $\eta_{z,t} = \eta_{f,t} = \eta_t$ we have that 
\begin{align}\label{eq:sec3:temp7}
    & \mathbb{E}[\bar{A}_{t+1}] =\sum_{i_z=1}^{N_z}  \mathbb{E}\left[\bar{A}_{i_z,t+1}\right] + \sum_{i_f=N_z+1}^{N}  \mathbb{E}\left[A_{i_f,t+1} \right]\nonumber \\
    &\leq  \sum_{i_z=1}^{N_z} \Big( \mathbb{E}  [\bar{A}_{i_z,t}] + \frac{|\delta_{t-1} - \delta_t|^2 D^2}{r^2} + \frac{2 |\delta_{t-1} - \delta_t| D^2}{r} \nonumber \\
    &\quad -  2 \eta_{t} \mathbb{E} \left\langle \nabla_{i_z} C_{i_z}(x_t)-\nabla_{i_z} \mathcal{C}_{i_z}(x_t^{*}),x_{i_z,t}- x_{i_z,t}^{*} \right\rangle \nonumber \\
    &\quad + \frac{\eta_{t}^2 d^2 U^2}{\delta_t^2}  + \frac{2\eta_{t} L_1 \delta_t D^2}{r} + \frac{2\eta_{t} L_0 \delta_t D }{r} \nonumber \\
    &\quad + 2 \eta_{t} L_1 \delta_t \sqrt{N_z} D \Big) \nonumber \\
    &\quad  + \sum_{i_f=N_z+1}^{N} \Big( \mathbb{E}[A_{i_f,t}] + \eta_{t}^2 L_0^2 + 2 \eta_{t} L_1 \delta_t \sqrt{N_z} D \nonumber \\
    &\quad  -  2 \eta_{t} \mathbb{E}\left\langle \nabla_{i_f} \mathcal{C}_{i_f}(x_t)- \nabla_{i_f} \mathcal{C}_{i_f}(x_t^{*}),x_{i_f,t}- x_{i_f}^{*} \right\rangle  \nonumber \\
    &\quad + \frac{2 \eta_{t} L_1 \delta_t D ^2 \sqrt{N_z}}{r}\Big) \nonumber \\
    &=  \mathbb{E}[\bar{A}_t] +  N_z \frac{|\delta_{t-1} - \delta_t|^2 D^2}{r^2} + \frac{2 N_z  |\delta_{t-1} - \delta_t| D^2}{r}
     \nonumber \\
    & \quad -  2 \eta_{t} \sum_{i=1}^N \mathbb{E}\left\langle \nabla_{i} C_{i}(x_t)-\nabla_{i} \mathcal{C}_{i}(x_t^{*}),x_{i,t}- x_{i,t}^{*} \right\rangle \nonumber \\
    &\quad + N_z \Big(\frac{ \eta_{t}^2 d^2 U^2}{\delta_t^2}  + \frac{2\eta_{t} L_1 \delta_t D^2}{r} + \frac{2\eta_{t} L_0 \delta_t D }{r} \nonumber \\
    &\quad \quad + 2 \eta_{t} L_1 \delta_t \sqrt{N_z} D \Big)\nonumber \\
    &\quad +(N-N_z) \Big(   \eta_{t}^2 L_0^2 + 2 \eta_{t} L_1 \delta_t \sqrt{N_z} D    + \frac{2 \eta_{t} L_1 \delta_t D ^2 \sqrt{N_z}}{r}       \Big)\nonumber \\
    &\leq  (1-2m \eta_t)\mathbb{E}[\bar{A}_t]  +\frac{2 N_z  |\delta_{t-1} - \delta_t| D^2}{r}    \nonumber \\
    & \quad + \frac{N_z |\delta_{t-1} - \delta_t|^2 D^2}{r^2}  + \frac{ \eta_{t}^2 N_z d^2 U^2}{\delta_t^2} +(N-N_z)\eta_{t}^2 L_0^2 \nonumber \\
    &\quad + 2 \eta_t \delta_t N L_1  (\frac{D^2 \sqrt{N_z} }{r} + \sqrt{N_z}D) + 2\eta_t \delta_t N_z L_0 \frac{D}{r},
\end{align}
where the last inequality follows from the strong monotonicity condition \eqref{eq:strong_monotone}. 
Recalling $\delta_t = \frac{N_z^{\frac{1}{6}}}{N^{\frac{1}{3}} t^{\frac{1}{3}}}$, it is easy to verify that 
\begin{align}\label{eq:sec3:temp7.1}
    |\delta_{t-1} - \delta_t|\leq \frac{N_z^{\frac{1}{6}}}{N^{\frac{1}{3}}  t^{\frac{4}{3}}},
\end{align}
for all $t\geq 1$. 
The inequality \eqref{eq:sec3:temp7.1} holds at $t=1$ since $|\delta_0-\delta_1| =\delta_1$.
Substituting $\eta_t = \frac{1}{mt}$, $\delta_t = \frac{N_z^{\frac{1}{6}}}{N^{\frac{1}{3}} t^{\frac{1}{3}}}$ and \eqref{eq:sec3:temp7.1} into the inequality \eqref{eq:sec3:temp7}, we have that
\begin{align}\label{eq:sec3:temp8}
    &\mathbb{E}[\bar{A}_{t+1}] \nonumber \\
    &\leq  (1-\frac{2}{t} )\mathbb{E}[\bar{A}_t] +  \frac{2 N_z^{\frac{7}{6}} D^2} {rN^{\frac{1}{3}}  t^{\frac{4}{3}}} + \frac{N_z^{\frac{4}{3}} D^2}{ r^2 N^{\frac{2}{3}}  t^{\frac{8}{3}}}
    + \frac{ d^2 U^2 N_z^{\frac{2}{3}} N^{\frac{2}{3}} }{ m^2 t^{\frac{4}{3}}  } \nonumber \\
    &\quad+ \frac{2  N_z^{\frac{1}{6}} N^{\frac{2}{3}} L_1 }{mt^{\frac{4}{3}} }( \frac{D^2 \sqrt{N_z}}{r} + \sqrt{N_z}D) +(N-N_z)\frac{L_0^2}{m^2 t^2} \nonumber \\
    &\quad + \frac{2N_z^{\frac{7}{6}} L_0 D }{mr N^{\frac{1}{3}} t^{\frac{4}{3}} } \nonumber \\
    &\leq  (1-\frac{2}{t} )\mathbb{E}[\bar{A}_t] +  \frac{2 N_z^{\frac{2}{3}} N^{\frac{1}{6}} D^2} { r t^{\frac{4}{3}}} + \frac{N_z^{\frac{2}{3}} D^2}{ r^2  t^{\frac{8}{3}}}
    + \frac{ d^2 U^2 N_z^{\frac{2}{3}} N^{\frac{2}{3}} }{ m^2 t^{\frac{4}{3}}  } \nonumber \\
    &\quad+ \frac{2  N_z^{\frac{1}{6}} N^{\frac{2}{3}}  }{mt^{\frac{4}{3}} }( \frac{L_1 D^2 \sqrt{N_z}}{r} + \sqrt{N_z} L_1D + \frac{L_0 D}{r}) \nonumber \\
    &\quad +(N-N_z)\frac{L_0^2}{m^2 t^2} \nonumber \\
    & \leq  (1-\frac{2}{t} )\mathbb{E}[\bar{A}_t] + \frac{N_z^{\frac{2}{3}} N^{\frac{2}{3}} }{ t^{\frac{4}{3}}} \Big( \frac{2D^2}{r} +\frac{D^2}{r^2} + \frac{ d^2 U^2}{m^2}\Big) \nonumber \\
    &\quad+ \frac{2  N_z^{\frac{1}{6}} N^{\frac{2}{3}}  }{mt^{\frac{4}{3}} }\Big( \frac{L_1 D^2 \sqrt{N_z} }{r} + \sqrt{N_z} L_1D + \frac{L_0 D}{r}\Big) \nonumber \\
    &\quad +(N-N_z)\frac{L_0^2}{m^2 t^2},
\end{align}
where the second inequality follows since $N_z \leq N$.

In the following, we use the standard mathematical inductive method to prove that the following inequality holds for all $ t\geq1$:
\begin{align}\label{eq:inductive:hybrid}
    \mathbb{E}[\bar{A}_t] \leq   \frac{N_z^{\frac{2}{3}} N^{\frac{2}{3}} }{ t^{\frac{1}{3}}} \Big( \frac{2D^2}{r} +\frac{D^2}{r^2} + \frac{ d^2 U^2}{m^2}\Big)  +(N-N_z)\frac{L_0^2}{m^2 t} \nonumber \\
    + \frac{2  N_z^{\frac{1}{6}} N^{\frac{2}{3}} }{mt^{\frac{1}{3}} }\Big( \frac{L_1 D^2 \sqrt{N_z}}{r} + \sqrt{N_z} L_1D + \frac{L_0 D}{r}\Big) +\frac{D^2}{t}.
\end{align}
Obviously, the inequality \eqref{eq:inductive:hybrid} holds at the initial time step $t=1$ since $\bar{A}_1 \leq D^2$.
In what follows, we show that, if the statement \eqref{eq:inductive:hybrid} holds at time step $t$, it must also hold for the next time step $t+1$. 
To do so, we combine the inequalities \eqref{eq:sec3:temp8} and \eqref{eq:inductive:hybrid}, yielding that
\begin{align}\label{eq:sec3:temp9}
    &\mathbb{E}[\bar{A}_{t+1}] \nonumber\\
    &\leq  (1-\frac{2}{t} ) \Big(  \frac{N_z^{\frac{2}{3}} N^{\frac{2}{3}} }{ t^{\frac{1}{3}}} \Big( \frac{2D^2}{r} +\frac{D^2}{r^2} + \frac{ d^2 U^2}{m^2}\Big)  +(N-N_z)\frac{L_0^2}{m^2 t}   \nonumber \\
    &\quad + \frac{2  N_z^{\frac{1}{6}} N^{\frac{2}{3}} }{mt^{\frac{1}{3}} }\Big( \frac{L_1 D^2 \sqrt{N_z}}{r} + \sqrt{N_z} L_1D + \frac{L_0 D}{r}\Big) +\frac{D^2}{t} \Big)  \nonumber \\
    &\quad +  \frac{N_z^{\frac{2}{3}} N^{\frac{2}{3}} }{ t^{\frac{4}{3}}} \Big( \frac{2D^2}{r} +\frac{D^2}{r^2} + \frac{ d^2 U^2}{m^2}\Big) +(N-N_z)\frac{L_0^2}{m^2 t^2} \nonumber \\
    &\quad + \frac{2  N_z^{\frac{1}{6}} N^{\frac{2}{3}} L_1 }{mt^{\frac{4}{3}} }\Big( \frac{L_1 D^2 \sqrt{N_z}}{r} + \sqrt{N_z} L_1D + \frac{L_0 D}{r}\Big)  \nonumber \\
    & = (1-\frac{1}{t} ) \Big(  \frac{N_z^{\frac{2}{3}} N^{\frac{2}{3}} }{ t^{\frac{1}{3}}} \Big( \frac{2D^2}{r} +\frac{D^2}{r^2} + \frac{ d^2 U^2}{m^2}\Big)  +(N-N_z)\frac{L_0^2}{m^2 t}   \nonumber \\
    &\quad + \frac{2  N_z^{\frac{1}{6}} N^{\frac{2}{3}} }{mt^{\frac{1}{3}} }\Big( \frac{L_1 D^2 \sqrt{N_z}}{r} + \sqrt{N_z} L_1D  + \frac{L_0 D}{r}\Big) \Big) \nonumber \\
    &\quad +(1-\frac{2}{t} ) \frac{D^2}{t}   \nonumber \\
    &\leq \frac{N_z^{\frac{2}{3}} N^{\frac{2}{3}} }{ (t+1)^{\frac{1}{3}}} \Big( \frac{2D^2}{r} +\frac{D^2}{r^2} + \frac{ d^2 U^2}{m^2}\Big)  +(N-N_z)\frac{L_0^2}{m^2 (t+1)} \nonumber \\ 
    &\quad + \frac{2  N_z^{\frac{1}{6}} N^{\frac{2}{3}} }{m (t+1)^{\frac{1}{3}} }\Big( \frac{L_1 D^2 \sqrt{N_z} }{r} + \sqrt{N_z} L_1D + \frac{L_0 D}{r}\Big) + \frac{D^2}{t+1},
\end{align}
where the last inequality follows since $(1-\frac{1}{t} ) \frac{1}{t^{\frac{1}{3}}  } \leq \frac{t}{t+1}\frac{1}{t^{\frac{1}{3}}  } \leq \frac{1}{(t+1)^{\frac{1}{3}}  }$ and $(1-\frac{2}{t})\frac{1}{t} \leq (1-\frac{1}{t})\frac{1}{t} \leq \frac{1}{t+1} $.
The inequality \eqref{eq:sec3:temp9} shows that the statement \eqref{eq:inductive:hybrid} also holds at $t+1$ given that it holds at $t$, and thus we can conclude that the statement \eqref{eq:inductive:hybrid} holds for all $ t\geq1$. 
Finally, we have
\begin{align*}
    &\mathbb{E}[A_T]  = \mathbb{E}[  \left\| x_T - x^{*} \right\|^2] \nonumber \\
    &= \mathbb{E}[  \left\| x_T - x_{T-1}^* + x_{T-1}^*- x^{*} \right\|^2] \nonumber \\
    &\leq \mathbb{E}[ 2 \bar{A}_T + 2  \left\|x_{T-1}^*- x^{*} \right\|^2] \nonumber \\
    &\leq \frac{2N_z^{\frac{2}{3}} N^{\frac{2}{3}} }{ T^{\frac{1}{3}}} \Big( \frac{2D^2}{r} +\frac{D^2}{r^2} + \frac{ d^2 U^2}{m^2}\Big)  +(N-N_z)\frac{2L_0^2}{m^2 T} \nonumber \\
    &\quad + \frac{4  N_z^{\frac{1}{6}} N^{\frac{2}{3}} L_1 }{mT^{\frac{1}{3}} }\Big( \frac{L_1 D^2 \sqrt{N_z}}{r} + \sqrt{N_z} L_1D + \frac{L_0 D}{r}\Big) +\frac{2D^2}{T} \nonumber \\
    &\quad + \frac{2 N_z^{\frac{1}{3}} D^2 }{N^{\frac{2}{3}} (T-1)^{\frac{2}{3}}  r^2} \nonumber \\
    &= \mathcal{O}\Big( N_z^{\frac{2}{3}} N^{\frac{2}{3}} T^{-\frac{1}{3}} +  (N-N_z)T^{-1}\Big),
\end{align*}
where the first inequality follows since $(a+b)^2 \leq 2 a^2 + 2b^2$.
The proof is complete.

\end{document}